\documentclass[times,sort&compress,3p]{elsarticle} 
\journal{arXiv}
\usepackage[labelfont=bf]{caption}
\usepackage{color}
\usepackage{amssymb,latexsym,bm,amsmath,amsthm,amsfonts}
\usepackage{eufrak, booktabs,comment,ulem}
\usepackage{multirow,mathrsfs}
\usepackage{epsfig}
\usepackage{enumitem}
\usepackage{lscape}
\usepackage{url}
\graphicspath{{./figures}}

\def\red{\color{red}}

\definecolor{purple}{rgb}{0.7,0.35,0.5}

\definecolor{green}{rgb}{0, 0.4, 0}
\include{bold}

\newcommand{\argmin}{\mathop{\rm argmin}}
\newcommand{\argmax}{\mathop{\rm argmax}}
\newcommand{\iid}{\stackrel{\rm iid}{\sim}}

\newcommand{\eqdef}{\stackrel{\rm def}{=}}
\newcommand{\convas}{\stackrel{\rm a.s.}{\rightarrow}}

\newcommand{\convd}{\stackrel{\rm d}{\rightarrow}}

\def\1{{\bm 1}}
\def\0{{\bm 0}}
\def\th{{\rm th}}

\def\cov{{\rm Cov}}

\def\gic{{\rm GIC}}

\def\I{\mathcal{I}}
\def\F{\mathcal{F}}

\def\th{{\rm th}}

\newcommand{\real}{\mathbb{R}}
\renewcommand{\tilde}{\widetilde}

\def\gic{{\rm GIC}}

\def\ica{{{\rm IC}_1}}
\def\icb{{{\rm IC}_2}}
\def\icc{{{\rm IC}_3}}
\def\pca{{{\rm PC}_1}}
\def\pcb{{{\rm PC}_2}}
\def\pcc{{{\rm PC}_3}}

\newtheorem{thm}{Theorem}

\newtheorem{lem}{Lemma}

\begin{document}
\title{Information Criterion-Based Rank Estimation Methods for Factor Analysis: \\ A Unified Selection Consistency Theorem and Numerical Comparison}
\begin{frontmatter}
\author[1]{Toshinari Morimoto\corref{cor1}}
\author[2]{Hung Hung}
\author[3]{Su-Yun Huang}
\address[1]{Department of Mathematics, National Taiwan University\\
No. 1, Sec. 4, Roosevelt Rd., Taipei 10617, Taiwan}
\address[2]{Institute of Health Data Analytics and Statistics, National Taiwan University\\
No.17 , Xu-Zhou Rd., Taipei 10055, Taiwan}
\address[3]{Institute of Statistical Science, Academia Sinica\\
No.128, Academia Road, Section 2, Nankang, Taipei 11529, Taiwan}
\cortext[cor1]{Corresponding author. Email: d09221002@ntu.edu.tw}
\date{}

\begin{abstract}
Over the years, numerous rank estimators for factor models have been proposed in the literature.
This article focuses on information criterion-based rank estimators and investigates their consistency in rank selection. 
The gap conditions serve as necessary and sufficient conditions for rank estimators to achieve selection consistency under the general assumptions of random matrix theory. 
We establish a unified theorem on selection consistency, presenting the gap conditions for information criterion-based rank estimators with a unified formulation.

To validate the theorem's assertion that rank selection consistency is solely determined by the gap conditions, we conduct extensive numerical simulations across various settings. 
Additionally, we undertake supplementary simulations to explore the strengths and limitations of information criterion-based estimators by comparing them with other types of rank estimators.
\end{abstract}


\end{frontmatter}

%
%

\section{Introduction}\label{sec:introduction}

\subsection{Background}\label{subsec:background}
Factor analysis (FA), initially introduced by \citet{spearman}, interprets observed data $\bm{x}_i \in \mathbb{R}^p$ as a composition of signal and noise components.
The signal component is represented by lower-dimensional data $\bm{f}_i$, expanded to $p$ dimensions through linear transformation.
In many studies and works of literature,, the factor model is commonly defined by the following equation:
\begin{equation}\label{eq:fa_model}
\bm{x}_i - \bm{\mu} = \bm{B} \bm{f}_i + \bm{\epsilon}_i, \quad i \in \{1,\ldots,n\},
\end{equation}
where $\bm{\mu} \in \mathbb{R}^p$ represents the mean vector of $\bm{x}_i$, $\bm{B}$ is a $p \times r$ matrix,
$\bm{f}_i \in \mathbb{R}^r$ and $\bm{\epsilon}_i \in \mathbb{R}^p$ are, respectively, mean-zero random vectors of factors and noises.
Moreover, assume the covariance matrix of $\bm{f}_i$ to be $\cov[\bm{f}_i] = \bm{I}_r$, and
$\cov[\bm{\epsilon}_i] = \bm{D} \eqdef \text{diag}(\sigma_1^2,\ldots,\sigma_p^2)$ represents the covariance matrix of $\bm{\epsilon}_i$. 
Also assume that $\bm{f}_i$ and $\bm{\epsilon}_i$ are stochastically independent.
The covariance matrix of $\bm{x}_i$ is $\bm{\Sigma} = \bm{B}\bm{B}^\top + \bm{D}$.

A significant challenge in FA involves estimating the ``{\it rank}'' of covariance matrix $\bm{\Sigma}$, which, in a broad sense, is often viewed as the number of factors. 
However, the rigorous definition of rank is given within the framework of random matrix theory (RMT), independent of the assumptions of factor models, in high-dimensional settings where $n, p \rightarrow \infty$.
This RMT-based definition of rank can also be applied as a generalization of the number of principal components in PCA.
In recent years, there has been the development of several rank estimators, roughly classified into three categories,
aiming to provide accurate estimations of the number of factors in FA and the number of principal components in PCA.
\begin{enumerate}\itemsep=0pt
\item Methods based on information criterion \cite{bai_ng2002, bai2018, hung2022}.
\item Methods based on upper bound estimation of the limiting spectral distribution of the sample covariance matrix \cite{dpa, bema, act}.
\item Methods based on the difference of successive eigenvalues of the sample covariance matrix \cite{on, ed, er_gr}.
\end{enumerate}

\subsection{Purpose and contributions}\label{subsec:purpose_and_contributions}

Our study primarily centers on estimators in Category 1, i.e., information criterion-based rank estimators.
Previous research by \citet{bai2018} laid the groundwork for gap conditions, which are the necessary and sufficient conditions to achieve the rank selection consistency, for AIC and BIC.
Building on this, \citet{hung2022} also derived the gap conditions for GIC under more generalized assumptions.

This article presents an extension of the work by \citet{bai_ng2002} by deriving the gap conditions for their PC$_1$, IC$_1$, PC$_2$, IC$_2$, PC$_3$ and IC$_3$ under the general RMT assumptions.
By combining these results with the findings of \citet{bai2018} and \citet{hung2022}, we establish a unified selection consistency theorem for information criterion-based estimators.
To demonstrate the validity of the theorem, we conduct numerical simulations.
Furthermore, we perform a comparative analysis of information criterion-based estimators against selected estimators from Categories 2 and 3,
aiming to explore the advantages of information criterion-based estimators as well as uncover any potential disadvantages.
The contributions of our study are as follows:
\begin{itemize}\itemsep=0pt
\item
A unified selection consistency theorem: 
We derive the gap conditions, the necessary and sufficient conditions for rank selection consistency, for $\pca, \ica, \pcb, \icb, \pcc$ and $\icc$ in \cite{bai_ng2002} under the general RMT assumptions.
We combine our results with prior research to present a unified selection consistency theorem for information criterion-based rank estimators.
\item
Validation of the theorem through simulations: 
We conduct numerical simulations to support our claim that the selection consistency of the information criterion-based estimators is determined by the gap conditions.
\item
Comparative analysis: 
We draw comparisons between information criterion-based estimators and some selected methods from Categories 2 and 3, to explore the unique features of the information criterion-based estimators as well as their potential limitations.
\end{itemize}

This article makes substantial contributions to the understanding and development of rank estimation methods in FA and PCA, with a specific focus on information criterion-based estimators.
Our work offers valuable insights and guidance for researchers and practitioners alike.

\subsection{Organization}
The remainder of this article is organized as follows.
Section~\ref{sec:rank_estimation_problem} gives an introduction to RMT used to analyze the asymptotic behavior of rank estimators.
We also illustrate a rigorous definition of rank grounded in the principles of RMT.
Section~\ref{sec:review_of_info_based_methods} provides a brief review of information criterion-based rank estimators.
Our main theoretical contribution, a unified selection consistency theorem, is presented in Section~\ref{sec:unified_theorem}.
This theorem establishes the gap conditions as necessary and sufficient conditions for rank selection consistency for information criteria, providing a nuanced understanding of their role in rank estimation.
Section~\ref{sec:numerical_studies} validates the unified selection consistency theorem through numerical simulations.
Additionally, we compare these information criterion-based estimators with estimators from Categories 2 and 3, underscoring their potential advantages and limitations.
Finally, Section~\ref{sec:conclude} summarizes the key findings and contributions of our study, suggesting potential directions for future research in the field of rank estimation.
All proofs are compiled in Appendix~\ref{sec:proof_of_theorem} for easy reference.
Appendix~\ref{sec:study3} provides additional simulation results, particularly focusing on the performance for BIC, PCs, and ICs in scenarios where $p\gg n$.
For readers' convenience, some selected rank estimation methods from Categories~2 and~3 are briefly reviewed in Appendix~\ref{sec:intro_category2_and_3}.

\subsection{Terms and Notations}\label{subsec:terms_notations}

We denote the observed data as $\bm{x}_1, \ldots, \bm{x}_n \in \mathbb{R}^p$, with the population mean and covariance matrix represented by $\bm{\mu} \eqdef {\rm E}[\bm{x}_i]$ and $\bm{\Sigma} \eqdef \cov[\bm{x}_i]$, respectively.
The population eigenvalues of $\bm{\Sigma}$, arranged in descending order, are denoted by $\{\lambda_j\}_{j=1}^p$.
The sample mean is denoted by $\bar{\bm{x}} \eqdef \frac{1}{n} \sum_{i=1}^n \bm{x}_i$, and the sample covariance matrix is denoted by $\bm{S}_n \eqdef \frac{1}{n-1} \sum_{i=1}^n (\bm{x}_i-\bar{\bm{x}}) (\bm{x}_i-\bar{\bm{x}})^\top$.
The sample eigenvalues of $\bm{S}_n$, arranged in descending order, are denoted by $\{\widehat{\lambda}_j\}_{j=1}^p$. 
The spectral decompositions of $\bm{\Sigma}$ and $\bm{S}_n$ are represented by $\bm{\Gamma} \bm{\Lambda} \bm{\Gamma}^\top$ and $\widehat{\bm{\Gamma}} \widehat{\bm{\Lambda}} \widehat{\bm{\Gamma}}^\top$, respectively. Here $\bm{\Gamma} = [\bm{\gamma}_1, \ldots, \bm{\gamma}_p]$, $\bm{\Lambda} = \text{diag}(\lambda_1, \ldots, \lambda_p)$, $\widehat{\bm{\Gamma}} = [\widehat{\bm{\gamma}}_1, \ldots, \widehat{\bm{\gamma}}_p]$, and $\widehat{\bm{\Lambda}} = \text{diag}(\widehat{\lambda}_1, \ldots, \widehat{\lambda}_p)$.
Although $\lambda_j$ and $\widehat{\lambda}_j$ should technically be indexed as $\lambda_j^{(p)}$ and $\widehat{\lambda}_j^{(n,p)}$ to reflect their dependence on $p$ and $(n,p)$, we simplify them to $\lambda_j$ and $\widehat{\lambda}_j$ to avoid overly complicated notations.
The empirical spectral distribution (ESD) of a matrix $\bm{A}$, as outlined in Section~\ref{subsec:rmt}, is denoted by $F^{\bm{A}}$. When a sequence of ESDs weakly converges to a distribution function, we refer to this distribution as a limiting spectral distribution (LSD).

%
%

\section{RMT fundaments and the concept of rank}\label{sec:rank_estimation_problem}

%
%
\subsection{Review of RMT}\label{subsec:rmt}
The review materials pertinent to the rank estimation problem are mainly taken from \citet{yao_zheng_bai_2015} and \citet{bai2018}.
Let $\bm{A}$ be a $p \times p$ symmetric positive semi-definite matrix and let $\{\alpha_j\}_{j=1}^p$ be its eigenvalues in descending order.
The empirical spectral distribution (ESD) of $\bm{A}$ is defined as
\begin{equation}\label{eq:def_ESD}
F^{\bm{A}}(t) \eqdef \frac{1}{p} \sum_{j=1}^p \I(\alpha_j \leq t).
\end{equation}
The following general assumptions are commonly imposed in the field of RMT.
\begin{itemize}\itemsep=0pt
\item[(C1)] $\bm{x}_i = \bm{\mu} + \bm{\Sigma}^{1/2} \bm{z}_i$, where the entries of $\bm{z}_i$ are i.i.d. random variables with mean $0$, variance $1$ and finite fourth moment.
\item[(C2)] $p/n\rightarrow c \in (0,\infty)$.
\item[(C3)] $F^{\bm{\Sigma}} \rightarrow H$ weakly as $p \rightarrow \infty$, where the support of the distribution $H$
is a bounded interval $[\underline{\lambda},\overline{\lambda}]$ with $\underline\lambda\ge 0$.
\end{itemize}
Under the general RMT assumptions (C1)-(C3), $F^{\bm{S}_n}(t)$, the ESD of $\bm{S}_n$, converges to the generalized Marchenko-Pastur (MP) distribution $\F_{c,H}$ indexed by $(c,H)$ with probability one. 
The distribution $\F_{c,H}$ is characterized by its Stieltjes transform $s_{c,H}(z)\eqdef\int\frac{1}{t-z} d{\F_{c, H}(t)}$, which satisfies the following equation:
\begin{eqnarray}
s_{c,H}(z) = \int \frac{1}{t\{1-c-czs_{c, H}(z)\}-z}dH(t), \label{eq:MPdist_characterization}
\end{eqnarray}
where $z \in \mathbb{C}^{+} \eqdef \{x+iy\mid x \in \mathbb{R}, y>0\}$.
Under high-dimensional settings where both $n$ and $p$ tend to infinity, $\widehat{\lambda}_j$ is not a consistent estimator of $\lambda_j$. 
Instead, we have the following properties:
\begin{description}\itemsep=0pt
\item[Case 1.] 
When $\lambda_j$ tends to infinity, we have
$\widehat{\lambda}_j/\lambda_j \convas 1$.
\item[Case 2.] 
When $\lambda_j$ converges to $\lambda_j^* \notin [\underline{\lambda}, \overline{\lambda}]$ with $\psi'(\lambda_j^*)>0$, we have
$\widehat{\lambda}_j \convas \psi(\lambda_j^*)$.
\item[Case 3.] 
When $\lambda_j$ converges to $\lambda_j^* \notin [\underline{\lambda}, \overline{\lambda}]$ with $\psi'(\lambda_j^*) \leq 0$, we have
$\begin{cases} \widehat{\lambda}_j \convas a, & {\rm for}~\lambda_j^* < \underline{\lambda}, \\  \widehat{\lambda}_j \convas b, & {\rm for}~\lambda_j^* > \overline{\lambda}. \end{cases}$
\item[Case 4.]
When $\lambda_j$ converges to $\lambda_j^* \in [\underline{\lambda}, \overline{\lambda}]$, we have $\widehat{\lambda}_j \convas u_{1-\alpha}$ where $\alpha \eqdef \lim_{p \to \infty} j/p$.
\end{description}
Here $\psi(\lambda) \eqdef \lambda \left\{ 1+c \int \frac{t}{\lambda-t} dH(t)\right\}$, 
$a$ and $b$ are the essential inifimum and the essential supremum of $\F_{c,H}$ respectively, 
and $u_{1-\alpha}$ is the $(1-\alpha)^{\rm th}$ quantile of $\F_{c,H}$.
The function $\psi(\lambda)$ is defined outside the support of $H$, that is, $(-\infty, \underline{\lambda}) \cup (\overline{\lambda}, \infty)$.
According to \citet[p.221]{yao_zheng_bai_2015}, there are two points at which $\psi'(\lambda)=0$.
These points are denoted as $\lambda_a$ and $\lambda_b$, with $\lambda_a < \lambda_b$.
The point $\lambda_a$ is located within $(-\infty, \underline{\lambda})$ and $\lambda_b$ is within $(\overline{\lambda}, \infty)$.
Moreover, $\psi(\lambda)$ is concave in the interval $(-\infty, \underline{\lambda})$, transitioning from $-\infty$ to $\psi(\lambda_a)$ and back to $-\infty$, forming a ``$\cap$" shape.
Conversely, $\psi(\lambda)$ is convex in the interval $(\overline{\lambda}, \infty)$, transitioning from $+\infty$ to $\psi(\lambda_b)$ and back to $+\infty$, forming a ``$\cup$" shape.
The essential supremum $b$ of $\F_{c,H}$ is represented by $\psi(\lambda_b)$, and similarly, the essential infimum $a$ is given by $\psi(\lambda_a)$. 
Particularly in the case when $j=o(p)$, it is important to note that if $\lambda_j$ converges to a finite number $\lambda_j^*$ with $\psi'(\lambda_j^*) > 0$, then it is identifiable. Conversely, if $\lambda_j$ converges to $\lambda_j^*$ but with $\psi'(\lambda_j^*) \leq 0$, then $\widehat{\lambda}_j$ converges to $b$, which is irrelevant to $\lambda_j$, rendering $\lambda_j$ non-identifiable.
In cases where $\lambda_j$ diverges to infinity, we have $\widehat{\lambda}_j/\lambda_j\convas 1$. Thus, diverging $\lambda_j$'s are considered identifiable.
These observations lead to defining the rank as the count of identifiable population eigenvalues among those for which $j = o(p)$, as elaborated in Section~\ref{subsec:definition_of_rank}.

%
\subsection{Definition of rank}\label{subsec:definition_of_rank}

We define the rank as the number of identifiable population eigenvalues among those where $j=o(p)$ as shown below:
\begin{equation}\label{eq:rmt_rank_definition}
r_0(p) \eqdef \sup\left\{1 \leq j \leq p \mid \psi'(\lambda_j) > 0 ~{\rm and} ~ j=o(p) \right\}.
\end{equation}
It is worth noting that those $\lambda_j$'s within the support of $H$ are excluded because $\psi'(\cdot)$ is well-defined only for $\lambda \notin [\underline{\lambda}, \overline{\lambda}]$. 
This rank definition, grounded in RMT, is independent of factor model assumptions and the specific distribution of $\{\bm{x}_i\}_{i=1}^n$. 
It is determined solely by $c$ and the asymptotic behavior of $\{\lambda_j\}_{j=1}^p$. 
While the rank $r_0(p)$ could oscillate or diverge as $p$ increases, we assume it remains constant for sufficiently large $p$ in the remainder of this article.

The rank defined above counts the number of eigenvalues that are sufficiently large to be distinct from the distribution $H$, the LSD of $\bm{\Sigma}$. 
This definition not only plays a role in generalizing the number of factors in FA but also can be applied to defining the number of principal components in PCA in high-dimensional settings.
To demonstrate that $r_0$ defined above is a reasonable generalization of the number of factors $r$ in suitable factor models, the relationship between $r_0$ and $r$ is further discussed in Section~\ref{subsec:rank_vs_number_of_factors} below.

%
%
\subsection{Relationship between the rank of covariance matrix and the number of factors}\label{subsec:rank_vs_number_of_factors}

Rank estimation originally arose from the need to estimate the number of factors, and the rank gains a mathematically rigorous definition within the framework of RMT, as highlighted in Section \ref{subsec:definition_of_rank}.
The RMT-based rank definition plays an important role in generalizing the notion of the number of factors and provides a reasonable characterization.
In the subsequent discussion, we argue that, under suitable conditions, the number of factors and the rank established through RMT coincide. 

Let us consider a factor model (\ref{eq:fa_model}) with a fixed number of factors $r$, and let $r_0$ represent the rank defined based on RMT as discussed in Section \ref{subsec:definition_of_rank}.
Assume that the variances of the noise components $\{\sigma_j^2\}_{j=1}^p$ are ordered in descending order, and that $\sigma_j^2 \leq \overline{\lambda}$ holds for all $j$.
Let $\{\ell_j\}_{j=1}^p$ represent the eigenvalues of the signal component $\bm{B}\bm{B}^\top$, arranged in descending order as well.
(Note that $\sigma_j^2$'s and $\ell_j$'s are indexed by both $j$ and $p$, similar to $\lambda_j$'s.)
We have $\ell_{r+1}, \ldots, \ell_{p} = 0$ as the rank of $\bm{B}\bm{B}^\top$ is at most~$r$.
Furthermore, according to \citet[Theorem A.43]{bai_silverstein2010}, the ESD of $\bm{D}$ converges to $H(t)$ by the subsequent inequality:
\begin{equation}
\sup_{t} \left\vert F^{\bm \Sigma}(t)-F^{\bm D} (t)\right\vert \leq \frac{1}{p} {\rm rank}(\bm{B}\bm{B}^\top) \leq \frac{r}{p} \rightarrow 0.
\end{equation}
Note that, as can be seen from the above equation, $H$ is determined solely by the asymptotic behavior of $\bm{D}$.
Thus, it is fair to say that $\bm{B}\bm{B}^\top$ (or $\{\ell_j\}_{j=1}^p$) has no impact on $H$.
According to the properties of eigenvalues of the sum of two symmetric matrices in \citet[p. 300]{yao_zheng_bai_2015}, which is also known as Weyl's inequality, the population eigenvalues $\{\lambda_j\}_{j=1}^p$, $\{\ell_j\}_{j=1}^p$, and $\{\sigma^2_j\}_{j=1}^p$ share the following relationships:
\begin{eqnarray}
\lambda_{p-i-j} &\geq& \ell_{p-i} + \sigma^2_{p-j}, \\
\lambda_{i+j-1} &\leq& \ell_i + \sigma^2_j,
\end{eqnarray}
where $i, j$ are non-negetive integers chosen so that the indices of $\lambda_{(\cdot)}, \ell_{(\cdot)}, \sigma^2_{(\cdot)}$ fall within the range of $\{1,\ldots,p\}$.
Based on these inequalities, we further deduce:
\begin{eqnarray}
\lambda_{r}~ &\geq& \ell_r + \sigma^2_p, \label{eq:lambda_r_vs_l_r_sigma2_p} \\
\lambda_{r+1} &\leq& \sigma_1^2. \label{eq:lambda_rplus1_vs_sigma2_1}
\end{eqnarray}
By Equation (\ref{eq:lambda_r_vs_l_r_sigma2_p}), it can be seen that when $\ell_r$ is substantial enough to ensure that $\lambda_{r} > \lambda_b$, we have $\psi'(\lambda_1), \ldots, \psi'(\lambda_r) > 0$.
Note that $\lambda_b$ is not influenced by $\{\ell_j\}_{j=1}^p$ as it is determined by $c$ and $H$; this allows an increase in $\ell_r$ to make $\lambda_r$ greater than $\lambda_b$ without affecting $\lambda_b$.
Additionally, considering Equation (\ref{eq:lambda_rplus1_vs_sigma2_1}) and the assumption that $\sigma_1^2 \leq \overline{\lambda}$, we observe that $\lambda_{r+1}$ is within the support of $H$ ($\psi'$ is not well-defined), and thus not identifiable.
Since only $\lambda_1, \ldots, \lambda_r$ satisfy $\psi'(\cdot) > 0$, therefore, $r = r_0$ holds.

In conclusion, the number of factors $r$ and the rank $r_0$ defined by RMT coincide if the following two conditions are met:
(1) The $r^\th$ eigenvalue of $\bm{B}\bm{B}^\top$, $\ell_r$, is sufficiently large to make $\lambda_r$, the $r^\th$ eigenvalue of $\bm\Sigma$, an identifiable population eigenvalue. 
(2) All noise variances stay below the upper bound of the support of $H$. 

\citet{wang_fan2017} also discussed the relationship between the factor model and the spiked covariance model, a model with several leading population eigenvalues significantly larger than the others. 
They showed that, under the pervasiveness assumption on $\bm B$, the factor model becomes a rank-$r$ spiked covariance model with diverging leading eigenvalues of order $O(p)$. 
Note that, by {\bf Case~1} discussed above, diverging eigenvalues are identifiable and can be separated from the support of $H$. 
Thus, the number of factors coincides with the number of identifiable eigenvalues. 
Unlike  \cite{wang_fan2017}, our discussions above regarding the relationship between the rank and the number of factors do not assume the pervasiveness assumption and do not rely on diverging leading eigenvalues.

%
%

\section{Review of information criterion-based estimators}\label{sec:review_of_info_based_methods}
%
%
\subsection{AIC, BIC and GIC}\label{subsec:intro_aic_bic_gic}

The development of rank estimators based on AIC, BIC, and GIC initiates by fitting the observed data to the rank-$r$ simple spiked working model:
\begin{equation}\label{eq:simple_spiked_model_under_gaussian}
\begin{cases}
\bm{x}_i \iid N(\bm{\mu}, \bm{\Sigma}_r), \quad \bm{\Sigma}_r \eqdef \bm{\Gamma}\bm{\Lambda}\bm{\Gamma}^\top + \sigma^2(\bm{I}_p - \bm{\Gamma}\bm{\Gamma}^\top), \\
{\rm where} ~~ \bm{\Gamma}^\top \bm{\Gamma} = \bm{I}_r, \quad \bm{\Lambda} \eqdef {\rm diag}(\lambda_1, \ldots, \lambda_r).
\end{cases}
\end{equation}
It is important to underscore that while the working model described in Equation (\ref{eq:simple_spiked_model_under_gaussian}) serves as a foundation for the derivation of the estimators, it may not necessarily coincide with the true underlying data distribution.
However, the validity of the estimators is not restricted to this particular working model assumption.
As detailed in Section \ref{sec:unified_theorem}, the consistency of the derived estimators soley hinges on the adherence to the general RMT assumptions (C1)-(C3), independent of the specific working model outlined in Equation~(\ref{eq:simple_spiked_model_under_gaussian}).

The selection criteria, detailed in Equations (\ref{eq:aic_estimator})-(\ref{eq:gic_estimator}), are constructed by integrating penalty terms into the estimated log-likelihood.
Subsequently, the model minimizing these criteria is selected, thus providing an estimation for the rank:
\begin{eqnarray}
\widehat r_{\rm AIC}&\eqdef&\argmin_{r\le q}\ln|\widehat{\bm{\Sigma}}_r|+\frac{2}{n}b_r ,\label{eq:aic_estimator} \\
\widehat r_{\rm BIC}&\eqdef&\argmin_{r\le q}\ln|\widehat{\bm{\Sigma}}_r|+\frac{\ln(n)}{n}b_r ,\label{eq:bic_estimator} \\
\widehat r_{\rm GIC}&\eqdef&\argmin_{r\le q}\ln|\widehat{\bm{\Sigma}}_r|+\frac{2}{n}\widehat b_r^{\gic} ,\label{eq:gic_estimator}
\end{eqnarray}
where $q=o(p)$,
$\widehat{\bm{\Sigma}}_r \eqdef \widehat{\bm{\Gamma}}\widehat{\bm{\Lambda}}\widehat{\bm{\Gamma}}^\top + \widehat{\sigma}_r^2(\bm{I}_p - \widehat{\bm{\Gamma}}\widehat{\bm{\Gamma}}^\top)$,
$\widehat{\sigma}_r^2 \eqdef \frac{1}{p-r} \sum_{j=r+1}^p \widehat{\lambda}_j$,
$b_r \eqdef {pr-r(r+2)/2}+r+1+p$ represents the number of free parameters in (\ref{eq:simple_spiked_model_under_gaussian}), and
\begin{equation}\label{eq:gic_penalty_term}
\widehat b_{r}^\gic \eqdef
\left\{
{r \choose 2} +
\sum_{j\le r}\sum_{\ell>r}
\frac{\widehat\lambda_{\ell}(\widehat\lambda_{j}-\widehat\sigma_r^{2})}{\widehat\sigma_r^{2}(\widehat\lambda_{j}-\widehat\lambda_{\ell})}
\right\} + r +
\frac{\frac{1}{p-r}\sum_{j>r}\widehat\lambda_{j}^2}{(\frac{1}{p-r}\sum_{j>r}\widehat\lambda_{j})^2}.
\end{equation}

%
%
\subsection{Bai and Ng's criteria}\label{subsec:intro_pc_ic}

The estimators proposed by \citet{bai_ng2002} involves fitting the observed data to a rank-$r$ factor model~(\ref{eq:fa_model}).
Instead of employing the maximum likelihood estimation, the authors estimated the factor loadings $\bm{B}$ by minimizing the Frobenius norm of the estimated noise components.
Although our notation differs from the original notation used in \cite{bai_ng2002}, when mapped to our notation in Equation (\ref{eq:fa_model}), the minimization of the factor loadings is depicted as in Equation (\ref{eq:estimation_of_B_in_rank_r_fa_model}).
\begin{equation}\label{eq:estimation_of_B_in_rank_r_fa_model}
\begin{cases}
\argmin_{\bm{B}}\frac{1}{p}\|\bm{S}_n-\bm{B}\bm{B}^\top\|_F^2 = \left[\widehat{\lambda}_1^{1/2}\widehat{\bm{\gamma}}_1,\ldots, \widehat{\lambda}_r^{1/2}\widehat{\bm{\gamma}}_r \right]
, \\
\min_{\bm{B}}\frac{1}{p}\|\bm{S}_n-\bm{B}\bm{B}^\top\|_F^2 = \frac{1}{p} \sum_{j=r+1}^p \widehat{\lambda}_j,
\end{cases}
\end{equation}
where $\|\cdot\|_F$ is the Frobenius norm.
If one assumes $\bm{\mu}=\bm{0}$ as in \cite{bai_ng2002}, the definition of the sample covariance matrix $\bm{S}_n$ can be reinterpreted as $\frac{1}{n}\sum_{i=1}^n \bm{x}_i\bm{x}_i^\top$.
The authors proposed the following estimators for $r$, specifically $\pca$, $\pcb$, and $\pcc$:
\begin{eqnarray}
\widehat r_\pca&\eqdef&\argmin_{r\le q}\left(\frac{1}{p}\sum_{j=r+1}^p \widehat\lambda_j\right) + \widehat \sigma_q^2r\left(\frac{n+p}{np}\right) \ln \left(\frac{np}{n+p}\right) \label{eq:pc1_estimator}, \\
\widehat r_\pcb&\eqdef&\argmin_{r\le q}\left(\frac{1}{p}\sum_{j=r+1}^p\widehat\lambda_j\right) + \widehat \sigma_q^2r  \left(\frac{n+p}{np}\right) \ln (p\wedge n) \label{eq:pc2_estimator}, \\
\widehat r_\pcc&\eqdef&\argmin_{r\le q}\left(\frac{1}{p}\sum_{j=r+1}^p\widehat\lambda_j\right) + \widehat \sigma_q^2r \frac{\ln (p\wedge n)}{p\wedge n} \label{eq:pc3_estimator},
\end{eqnarray}
where $\widehat{\sigma}_q^2 \eqdef \widehat{\sigma}_{r \leftarrow q}^2$ as defined below (\ref{eq:gic_estimator}) is a consistent estimator of $\sigma^2 \eqdef \frac{1}{np}\sum_{i=1}^n E\left[\bm{\epsilon}_i^\top \bm{\epsilon}_i\right]$.
To avoid estimating $\sigma^2$, the authors further introduced the criteria $\ica$, $\icb$, and $\icc$:
\begin{eqnarray}
\widehat r_\ica&\eqdef&\argmin_{r\le q} \,\ln\left(\frac{1}{p}\sum_{j=r+1}^p \widehat\lambda_j\right) + r \left(\frac{n+p}{np}\right) \ln \left(\frac{np}{n+p}\right) \label{eq:ic1_estimator}, \\
\widehat r_\icb&\eqdef&\argmin_{r\le q} \, \ln\left(\frac{1}{p}\sum_{j=r+1}^p \widehat\lambda_j\right) + r \left(\frac{n+p}{np}\right) \ln (p\wedge n) \label{eq:ic2_estimator}, \\
\widehat r_\icc&\eqdef&\argmin_{r\le q} \, \ln\left(\frac{1}{p}\sum_{j=r+1}^p \widehat\lambda_j\right) + r \frac{\ln (p\wedge n)}{p\wedge n} \label{eq:ic3_estimator}.
\end{eqnarray}

\section{Main result: a unified selection consistency theorem}\label{sec:unified_theorem}
In this section, we establish a unified selection consistency theorem for the information-criterion based rank estimators. 
Here the term ``{\it selection consistency}" refers to each estimator $\widehat{r}$ converging to $r_0$ with probability 1, i.e., $\widehat{r} \convas r_0$. 
%
%
%
%
\begin{thm}[Selection consistency of information criterion-based estimators]\label{thm:info_gap_conditions}
Let $\psi_{r_0} \eqdef \psi(\lambda_{r_0})$ and $\mu_H \eqdef \int t dH(t)$.
Assume that the rank $r_0$, defined in~(\ref{eq:rmt_rank_definition}), is a fixed positive integer.
Under the general RMT assumptions (C1)-(C3), the information criterion-based estimators achieve selection consistency if and only if
\begin{equation}\label{eq:unified_gap_conditions}
\begin{cases}
\lim_{n,p\to\infty} L_{\rm m}(\psi_{r_0}/\mu_H)> 0 ,\\
\lim_{n,p\to\infty} L_{\rm m}(b/\mu_H) < 0 , 
\end{cases}
\end{equation}
where $L_{\rm m}(u) \eqdef g(u)-\beta_{\rm m}(u)$, $g(u) \eqdef u-1-\ln u$, and $\beta_{\rm m}(u)$ is a method-dependent function defined as follows:
\begin{itemize}\itemsep=0pt
\item
AIC corresponds to $\beta_{\rm AIC}(u) \eqdef 2c$.
\item
BIC corresponds to $\beta_{\rm BIC}(u) \eqdef c\ln n$.
\item
GIC corresponds to $\beta_{\rm GIC}(u) \eqdef 2\kappa(u)$,
where $\kappa(u) \eqdef c(u-1) \int \frac{t/\mu_H}{u-t/\mu_H} d\mathcal{F}_{c,H}(t)$.
\item
$\pca$ and $\ica$ correspond to $\beta_1(u) \eqdef g\left(\left(1+c\right)\ln(\frac{p}{1+c})\right)$.
\item
$\pcb$ and $\icb$ correspond to $\beta_2(u) \eqdef g\left(\left(1+c\right)\ln (n\wedge p)\right)$.
\item
$\pcc$ and $\icc$ correspond to $\beta_3(u) \eqdef g\left((1 \vee c) \ln (n\wedge p)\right)$.
\end{itemize}
\end{thm}
We refer to the first requirement in Equation (\ref{eq:unified_gap_conditions}) as ``{\it the first gap condition}'', and the second requirement as ``{\it the second gap condition}''.
These conditions delineate the minimum required signal $\psi_{r_0}$ (or $\lambda_{r_0}$) and the maximum acceptable noise $b$ for ensuring the selection consistency of each underlying estimator.
The gap conditions for AIC, BIC and GIC were derived and proven by \citet{bai2018} and \citet{hung2022}.
In Theorem~\ref{thm:info_gap_conditions}, we extend these results to include the gap conditions for $\pca$, $\ica$, $\pcb$, $\icb$, $\pcc$ and $\icc$, and present a unified selection consistency theorem.
This theorem reveals that the gap conditions for the information criterion-based estimators can be encapsulated within a unified mathematical formulation, with variations only in the method-dependent function $\beta_{\rm m}(\cdot)$.
The function $g(\cdot)$ quantifies the increase in log-likelihood for a simple spiked working model under the Gaussian working assumption (\ref{eq:simple_spiked_model_under_gaussian}) as the model rank increases by one unit.
The function $\beta_{\rm m}(\cdot)$ originates from the penalty term used in each estimator's derivation.
It should be emphasized that, for a given $(n,p,c)$, only $\beta_{\rm GIC}(u)$ depends on $u$, whereas the function $\beta_{\rm m}(\cdot)$ for other estimators remains constant.

It is noteworthy that despite the differing methodologies used to derive $\pca$, $\ica$, $\pcb$, $\icb$, $\pcc$, and $\icc$, compared to AIC, BIC and GIC, the gap conditions for all these estimators can be expressed through a unified mathematical formulation.
Remarkably, the gap conditions for $\pca$ and $\ica$ (as well as for $\pcb$ and $\icb$; $\pcc$ and $\icc$) are shown to be identical, indicating that $\pca$ and $\ica$ are likely to exhibit asymptotically similar properties.
A further consequence of this theorem is that the asymptotic behavior of these information criterion-based estimators is mainly governed by the function $\beta_{\rm m}(\cdot)$.
As $\beta_{\rm m}(\cdot)$ increases, the minimum required signal $\lambda_{r_0}$ also rises in the first gap condition, while permitting a higher level of noise $b$ in the second gap condition.
Hence, $\beta_{\rm m}(\cdot)$ plays a crucial role in determining both the required signal strength and the acceptable noise level.
Additional insights can be gained by examining the asymptotic behavior of $\beta_{\rm m}(\cdot)$.
The following summary outlines the fundamental aspects concerning the asymptotic behavior of these information criterion-based estimators:
\begin{itemize}\itemsep=0pt
\item
The following relationships hold: $\beta_{\rm AIC} < \beta_{\rm GIC} < \beta_{\rm BIC}$ and $\left(\beta_{\rm BIC} \vee \beta_{3}\right) <\beta_{1}<\beta_{2}$.
\item
The functions $\beta_{\rm BIC}$, $\beta_1$, $\beta_2$ and $\beta_3$ diverge as $n, p$ go to infinity.
Therefore, the selection consistency for BIC, $\pca/\ica$, $\pcb/\icb$, and $\pcc/\icc$ requires a diverging signal eigenvalue $\lambda_{r_0}$ in the first gap condition.
On the other hand, the second gap condition is automatically satisfied.
\item
The above observation for BIC, $\pca/\ica$, $\pcb/\icb$, and $\pcc/\icc$ ensures their robustness against a large noise $b$ but it also means that they might overlook weak signal eigenvalues.
This reveals a trade-off: researchers applying these estimators must weigh the sensitivity to signal eigenvalues against the robustness to noise.
\item
Among the functions $\beta_1$, $\beta_2$, and $\beta_3$, the relationship $\beta_3 < \beta_1 < \beta_2$ holds.
Additionally, $\beta_{\rm BIC} < \beta_1$ holds.
Moreover, as $c$ becomes sufficiently small, we have $\beta_1\approx \beta_2 \approx  \beta_3$.
When $c$ is either sufficiently large or sufficiently small, we have $\beta_2\approx \beta_1$.
\item
The relationship between $\beta_{\rm BIC}$ and $\beta_3$ varies depending on the scenarios of $c$.
When $c$ is small, $\beta_{\rm BIC} < \beta_3$ holds.
Conversely, for $c > 1$, we have $\beta_3 = g(\beta_{\rm BIC}) < \beta_{\rm BIC}$.
\end{itemize}


%
%

\section{Numerical studies}\label{sec:numerical_studies}

%
%
\subsection{Study I: verification of the gap conditions as necessary and sufficient conditions for selection consistency}\label{subsec:study1}

In Section \ref{subsec:study1}, we aim to demonstrate the validity of our paper's main theoretical result, the unified selection consistency theorem, through numerical experiments.
We generate samples $\bm{x}_1, \ldots, \bm{x}_n$ based on the various settings described in Section \ref{subsubsec:study1_settings}.
As discussed in Section \ref{subsec:gap_conditions}, for each specific setting, the satisfaction of the gap conditions for each estimator is determined.
Subsequently, in Section \ref{subsubsec:study1_results}, we examine if the ability of an estimator to perform accurate rank estimation under a specific setting hinges on these gap conditions.

\subsubsection{Experimental settings} \label{subsubsec:study1_settings}

We consider the following three LSDs for $\bm{\Sigma}$: $H_1$, $H_2$, and $H_3$, all of which have their support within $[0,1]$ and satisfy the requirement of the RMT assumption (C3).
Strictly speaking, the support of $H_2$ is not $[0,1]$; however, each point of mass of $H_2$ appears across the entirety of the interval $[0,1]$.
Therefore, for convenience, we choose to consider it as having support over $[0,1]$.
(The support of the distribution is defined by ${\rm supp}(H) \eqdef \left\{\lambda \mid {\rm Pr}(X \in (\lambda-\epsilon,\lambda+\epsilon)>0 ~ \forall \epsilon > 0)\right\}$, where $X \sim H$.)

\begin{itemize}\itemsep=0pt
\item $H_1$: A Beta distribution characterized by parameters $(\alpha,\beta)= (3,3)$.
\item $H_2$: A discrete distribution created by $\min(\frac{\Lambda+1}{50},1)$, where $\Lambda \sim {\rm Poisson}(24)$.
\item $H_3$: A truncated exponential distribution generated by $\min(\Lambda,1)$, with $\Lambda$ following an exponential distribution and ${\rm E}[\Lambda]=0.63$.
\end{itemize}
Note that $H_1$ is a continuous distribution, $H_2$ is a discrete distribution, and $H_3$ is a mixture of a continuous distribution and a discrete distribution with a point mass of one at $1$.
For the constant $c$, we set its values to 0.25, 0.5, 1, and 1.5, corresponding to $(n,p)$ values of $(1000,250)$, $(1000,500)$, $(600,600)$, and $(400,600)$, respectively.
We establish $\lambda_{r_0}$ as 2, 3, 4, and 5, and assign $r_0$ to 5.
The simulation procedure is as follows:
\begin{description}
\item[Step 1.]
Generate $\bm{\Lambda} \eqdef {\rm diag}(\lambda_1, \ldots, \lambda_p)$. 
Note that $\lambda_{r_0}$ is given as a condition in each individual setting.
The leading eigenvalues $\lambda_1, \ldots, \lambda_{r_0-1}$ are generated as $U_i + \lambda_{r_0}$, where $U_1, \ldots, U_{r_0-1}$ are i.i.d. random variables from a uniform distribution with support $(0,1)$.
The subsequent eigenvalues $\lambda_{r_0+1}, \ldots, \lambda_{p}$ are sampled from the given $H$ ($H_1$, $H_2$ or $H_3$).
\item[Step 2.]
Let $\bm{\Sigma} = \bm{\Gamma} \bm{\Lambda} \bm{\Gamma}^\top$, where $\bm{\Gamma}$ is drawn uniformly from the set of orthogonal $p \times p$ matrices.
\item[Step 3.]
Generate $\bm{x}_1, \ldots, \bm{x}_{n} \iid N(\bm{0}, \bm{\Sigma})$.
\item[Step 4.]
Estimate the rank using each rank estimation method. The upper bound of the search range for $r$ during rank search is set to  $q \eqdef \lfloor2\sqrt{\min\{n,p\}}\rfloor$.
\item[Step 5.]
Repeat the procedure above for $T=500$ times for each individual setting.
\end{description}

\subsubsection{Gap conditions}\label{subsec:gap_conditions}
Let us examine the satisfaction of the gap conditions for each estimator with respect to the settings described in the aforementioned Section \ref{subsubsec:study1_settings}.
As stated in Theorem \ref{thm:info_gap_conditions}, the gap conditions for each estimator under each setting are represented by Equation (\ref{eq:unified_gap_conditions}).
While noting that the function $\beta_{\rm m}(\cdot)$ is a constant function for methods other than GIC, gap conditions in Equation~(\ref{eq:unified_gap_conditions}) can be rewritten as follows:
($\lim$ is omitted)
\begin{equation}\label{eq:unified_gap_conditions_variant}
\begin{cases}
g\left(\frac{b}{\mu_H}\right)< \beta_{\rm m} < g\left(\frac{\psi_{r_0}}{\mu_H}\right) & {\rm for~AIC, BIC, PCs, ICs}, \\
g\left(\frac{b}{\mu_H}\right)< 2\kappa\left(\frac{b}{\mu_H}\right), \quad  2\kappa\left(\frac{\psi_{r_0}}{\mu_H}\right) < g\left(\frac{\psi_{r_0}}{\mu_H}\right)~ &{\rm for~GIC} .
\end{cases}
\end{equation}
We aim to provide a visual representation to better understand the satisfaction of gap conditions based on Equation (\ref{eq:unified_gap_conditions_variant}).
Given specific $c$ and $H$, the variables $b$ and $\mu_H$ are treated as constants, while $\psi_{r_0} \eqdef \psi(\lambda_{r_0})$ is expressed as a function of $\lambda_{r_0}$.
Consequently, with $c$ and $H$ predefined, we plot the curve $y=g(\frac{\psi_{r_0}}{\mu_H})$ and the horizontal line $y=g(\frac{b}{\mu_H})$ by designating $\lambda_{r_0}$ as the $x$-axis.
Note that when the $\beta_{\rm m}(\cdot)$ function includes $(n,p)$, the finite values of $(n,p)$ given in each individual setting are used in the following visual representation and also in the assessment of satisfaction of the gap conditions.

\begin{figure}[h!]
\includegraphics[width=12cm]{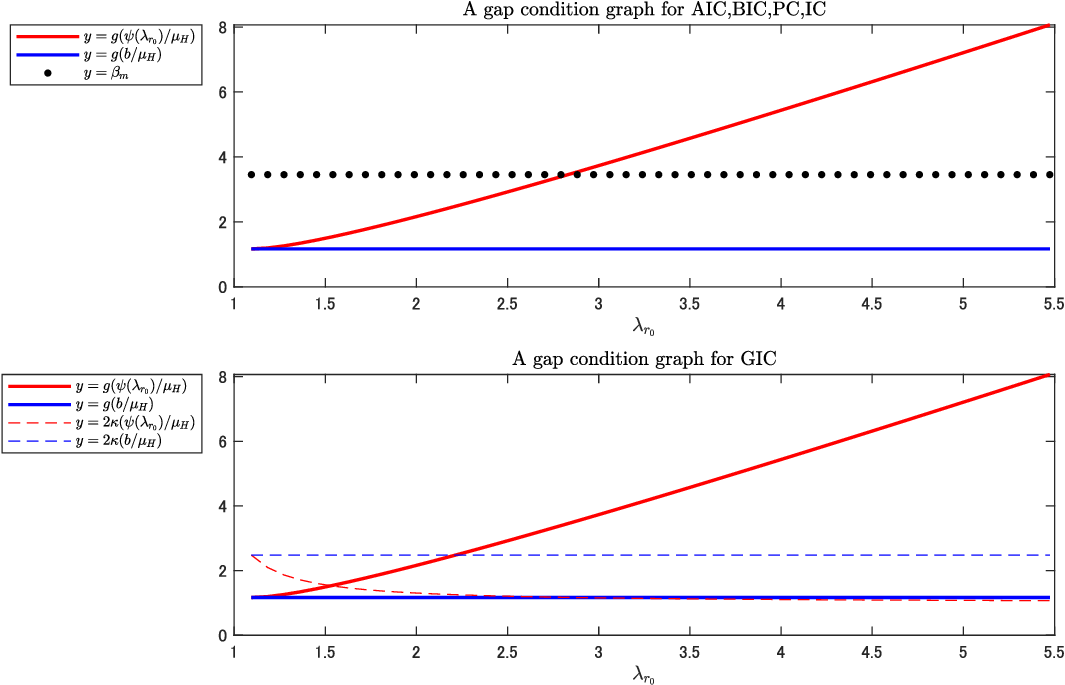}
\centering
\caption{
Graphs for the visual assessment of the satisfaction of the gap conditions as expressed in Equation (\ref{eq:unified_gap_conditions_variant}).
The upper graph is for methods of AIC, BIC, PCs, and ICs, whereas the lower one is for GIC.
In both, the red solid curve shows $y=g(\psi_{r_0}/\mu_H)$ and the blue solid horizontal line depicts $y=g(b/\mu_H)$.
For the upper graph, the gap conditions are satisfied where the plots for $y=\beta_{\rm m}$ are sandwiched between the red solid curve and blue solid line.
In the lower graph, the red dotted curve and blue dotted line represent $2\kappa(\psi_{r_0}/\mu_H)$ and $2\kappa(b/\mu_H)$, respectively.
Here, the gap conditions are satisfied when the red dotted curve is beneath the red solid curve, and the blue dotted line is above the blue solid line.
}
\label{fig:how_to_see_gap_condition_graphs}
\end{figure}

Figure~\ref{fig:how_to_see_gap_condition_graphs} is designed to visually assess whether the gap conditions are satisfied, comprising two graphs.
Both graphs contain the method-dependent $\beta_{\rm m}$ and two functions, $y=g\left(\frac{\psi_{r_0}}{\mu_H}\right)$ plotted with a red solid curve and $y=g\left(\frac{b}{\mu_H}\right)$ plotted with a blue solid horizontal line.
The upper graph is for methods of AIC, BIC, $\pca/\ica$, $\pcb/\icb$, and $\pcc/\icc$, where $\beta_{\rm m}$ is a constant function. 
The gap conditions are satisfied at values of $\lambda_{r_0}$ where $y=\beta_{\rm m}$ lies between the red solid curve and the blue solid horizontal line.
A failure to fall within this range indicates a violation of the gap conditions. 
Conversely, the lower graph clarifies the gap conditions for GIC.
Along with the red solid curve and the blue horizontal line, this graph introduces a red dotted curve and a blue dotted line, representing $2\kappa\left(\frac{\psi_{r_0}}{\mu_H}\right)$ and $2\kappa\left(\frac{b}{\mu_H}\right)$, respectively.
The gap conditions for GIC are satisfied when the red dotted curve is below the red solid curve as well as the blue dotted line is above the solid blue line.

\begin{figure}[h!]
\includegraphics[width=12cm, height=13cm]{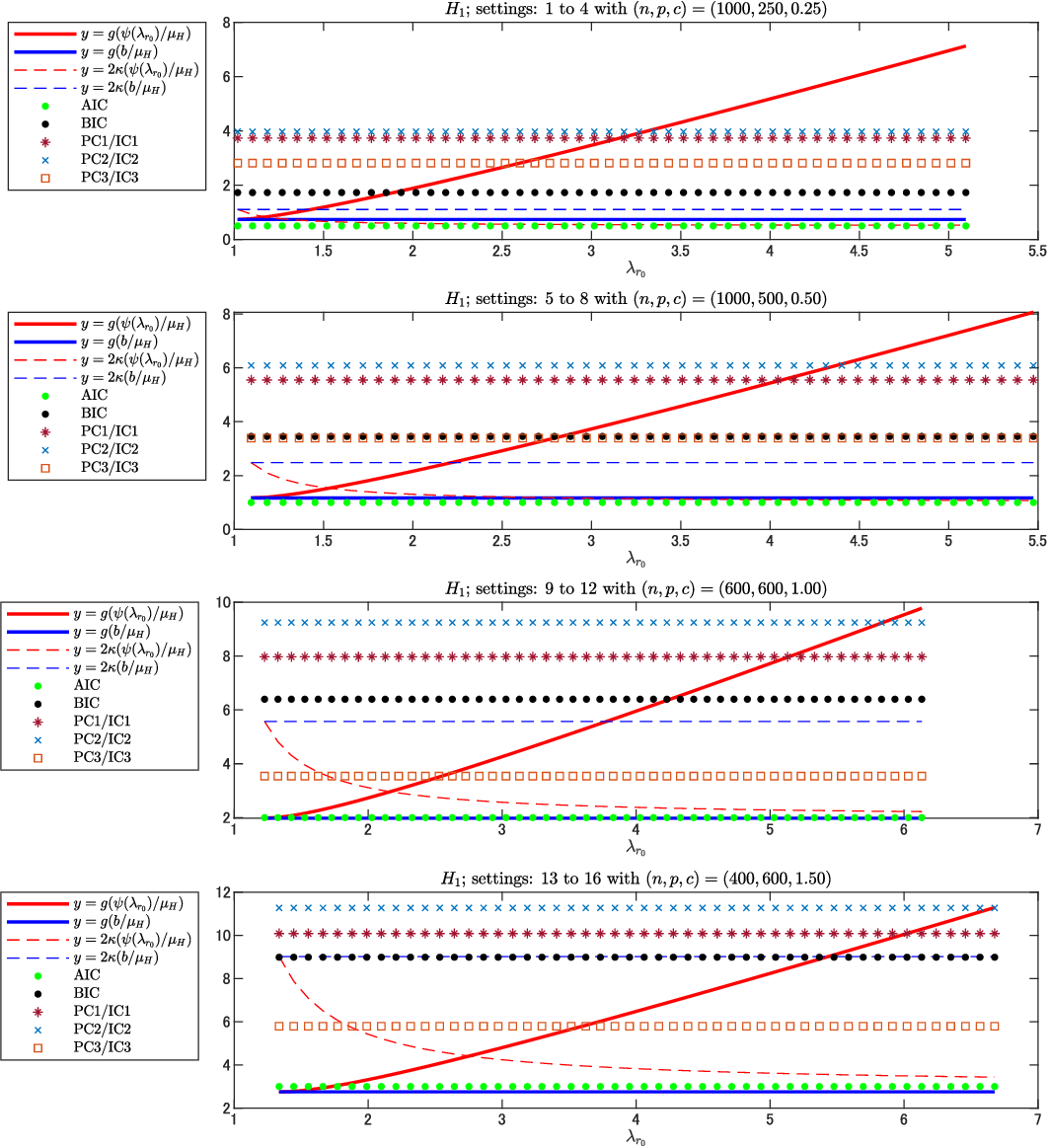}
\centering
\caption{$H_1$'s gap condition graphs:
The quartet of graphs corresponds to settings 1-4, 5-8, 9-12, and 13-16, respectively.
For information on which combinations of $n$, $p$, $\lambda_{r_0}$ correspond to each setting number, please refer to Table \ref{table:gap_conditions_H1}.
Detailed guidance on interpreting these graphs is elaborated at the beginning of Section \ref{subsec:gap_conditions}.
The plots represent the value of $\beta_{\rm m}$ for each estimator.
For AIC, BIC, PCs, and ICs, the gap conditions are satisfied within the range of $\lambda_{r_0}$ where the estimator's plots are enclosed between the red solid curve and the blue solid line.
For GIC, the gap conditions are satisfied within the range of $\lambda_{r_0}$ where the red dashed curve is beneath the red solid curve and the blue dashed line is above the blue solid line.
}
\label{fig:graph_of_gap_conditions_H1}
\end{figure}

\begin{figure}[h!]
\includegraphics[width=12cm, height=13cm]{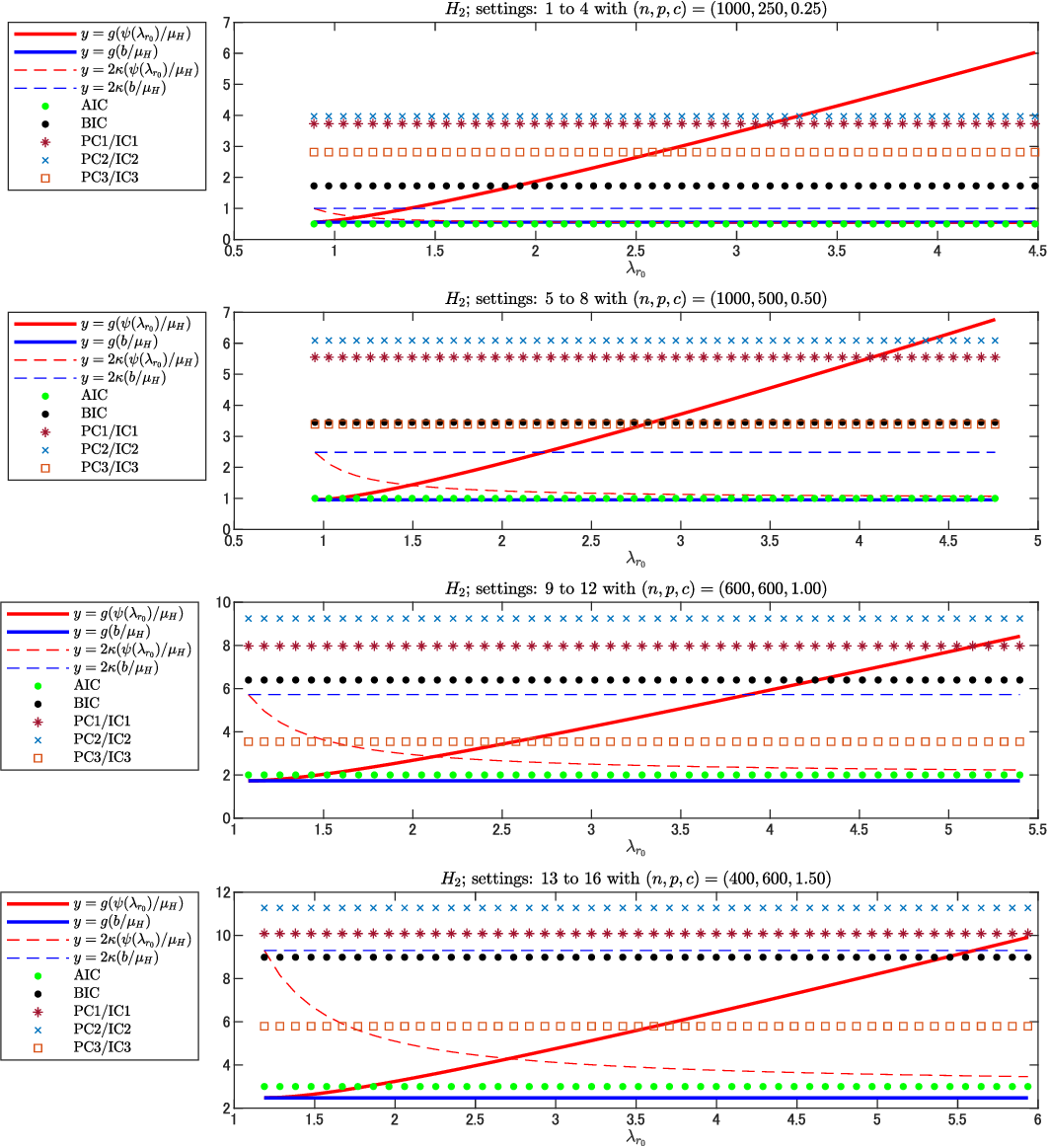}
\centering
\caption{
$H_2$'s gap condition graphs:
The quartet of graphs corresponds to settings 1-4, 5-8, 9-12, and 13-16, respectively.
For information on which combinations of $n$, $p$, $\lambda_{r_0}$ correspond to each setting number, please refer to Table \ref{table:gap_conditions_H2}.
Detailed guidance on interpreting these graphs is elaborated at the beginning of Section \ref{subsec:gap_conditions}.
The plots represent the value of $\beta_{\rm m}$ for each estimator.
For AIC, BIC, PCs, and ICs, the gap conditions are satisfied within the range of $\lambda_{r_0}$ where the estimator's plots are enclosed between the red solid curve and the blue solid line.
For GIC, the gap conditions are satisfied within the range of $\lambda_{r_0}$ where the red dashed curve is beneath the red solid curve and the blue dashed line is above the blue solid line.
}
\label{fig:graph_of_gap_conditions_H2}
\end{figure}

\begin{figure}[h!]
\includegraphics[width=12cm, height=13cm]{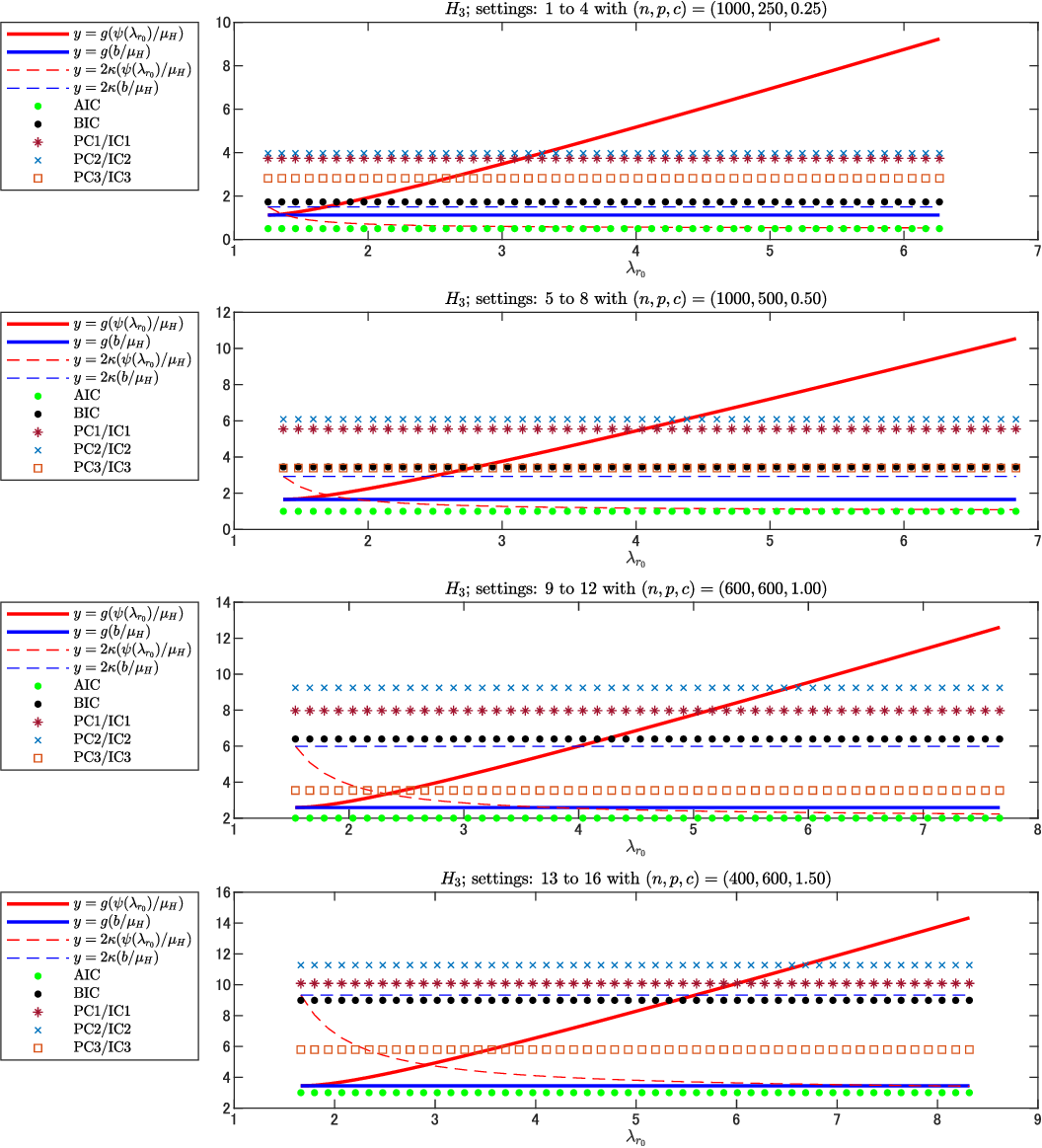}
\centering
\caption{
$H_3$'s gap condition graphs:
The quartet of graphs corresponds to settings 1-4, 5-8, 9-12, and 13-16, respectively.
For information on which combinations of $n$, $p$, $\lambda_{r_0}$ correspond to each setting number, please refer to Table \ref{table:gap_conditions_H3}.
Detailed guidance on interpreting these graphs is elaborated at the beginning of Section \ref{subsec:gap_conditions}.
The plots represent the value of $\beta_{\rm m}$ for each estimator.
For AIC, BIC, PCs, and ICs, the gap conditions are satisfied within the range of $\lambda_{r_0}$ where the estimator's plots are enclosed between the red solid curve and the blue solid line.
For GIC, the gap conditions are satisfied within the range of $\lambda_{r_0}$ where the red dashed curve is beneath the red solid curve and the blue dashed line is above the blue solid line.
}
\label{fig:graph_of_gap_conditions_H3}
\end{figure}

Figures \ref{fig:graph_of_gap_conditions_H1}-\ref{fig:graph_of_gap_conditions_H3} display graphs of the gap conditions for each simulation setting described in Section \ref{subsubsec:study1_settings}, drawn using the same methodology as in Figure \ref{fig:how_to_see_gap_condition_graphs}.
However, for comparative purposes, plots corresponding to all information criteria have been consolidated and displayed within a single graph.
Specifically, Figures \ref{fig:graph_of_gap_conditions_H1}-\ref{fig:graph_of_gap_conditions_H3} correspond to $H_1, H_2,$ and $H_3$, respectively.
Each of these figures comprises four individual graphs, each corresponding to a different combination of $(n,p)$.
The setting numbers listed in the captions of each graph (1-4, 5-8, 9-12,13-16) represent distinct combinations of $n, p,$ and $\lambda_{r_0}$ (refer to Tables \ref{table:gap_conditions_H1}-\ref{table:gap_conditions_H3} for details).

Tables~\ref{table:gap_conditions_H1}-\ref{table:gap_conditions_H3} present the values pertinent to the gap conditions as shown in Equation (\ref{eq:unified_gap_conditions_variant}).
Specifically, these tables include values such as $\psi_{r_0}$, $b$, $\mu_H$, along with $\beta_{\rm m}(\frac{\psi_{r_0}}{\mu_H})$ and $\beta_{\rm m}(\frac{b}{\mu_H})$ for each estimator denoted by $m$.
It is worth noting that for AIC, BIC, PCs, and ICs, $\beta_{\rm m}(\cdot)$ remains a constant function.
Corresponding values for $\beta_{\rm m}$ are articulated in columns labeled AIC, BIC, $\pca\ica$, $\pcb\icb$, and $\pcc\icc$.
For GIC, the values of $\beta_{\rm m}(\frac{\psi_{r_0}}{\mu_H})$ and $\beta_{\rm m}(\frac{b}{\mu_H})$ are provided in columns designated as $2\kappa(\psi_{r_0}/\mu_H)$ and $2\kappa(b/\mu_H)$, respectively.
The initial column delineates the specific setting number associated with a different combination of $n, p, \lambda_{r_0}$ under the given $H$.
Information summarizing whether each estimator satisfies the gap conditions in each setting is displayed along with the simulation results in Tables \ref{table:unified_theorem_simulation_H1}-\ref{table:unified_theorem_simulation_H3}.

\begin{table}[h!]
\scalebox{0.70}{
\begin{tabular}{|c|c|c|c|c|c|c|c|c|c|c|c|c|c|c|c|c|}
\hline
 & $n$  & $p$ & $c$  & $\lambda_{r_0}$ & $\psi_{r_0}$ & $b$   & $\mu_H$ & $g(\frac{\psi_{r_0}}{\mu_H})$ & $g(\frac{b}{\mu_H})$ & AIC & BIC   & $2\kappa(\frac{\psi_{r_0}}{\mu_H})$ & $2\kappa(\frac{b}{\mu_H})$ & $\pca\ica$  & $\pcb\icb$  & $\pcc\icc$ \\ \hline
1       & 1000 & 250 & 0.25 & 2               & 2.178        & 1.371 & 0.5     & 1.885                         & 0.734                & 0.5 & 1.727 & 0.59                                & 1.102                      & 3.732  & 3.97   & 2.813 \\ \hline
2       & 1000 & 250 & 0.25 & 3               & 3.155        & 1.371 & 0.5     & 3.472                         & 0.734                & 0.5 & 1.727 & 0.547                               & 1.101                      & 3.732  & 3.97   & 2.813 \\ \hline
3       & 1000 & 250 & 0.25 & 4               & 4.146        & 1.371 & 0.5     & 5.179                         & 0.734                & 0.5 & 1.727 & 0.531                               & 1.102                      & 3.732  & 3.97   & 2.813 \\ \hline
4       & 1000 & 250 & 0.25 & 5               & 5.141        & 1.371 & 0.5     & 6.949                         & 0.733                & 0.5 & 1.727 & 0.523                               & 1.101                      & 3.732  & 3.97   & 2.813 \\ \hline
5       & 1000 & 500 & 0.5  & 2               & 2.355        & 1.694 & 0.5     & 2.159                         & 1.167                & 1   & 3.454 & 1.306                               & 2.478                      & 5.549  & 6.09   & 3.388 \\ \hline
6       & 1000 & 500 & 0.5  & 3               & 3.31         & 1.694 & 0.5     & 3.73                          & 1.167                & 1   & 3.454 & 1.157                               & 2.478                      & 5.549  & 6.09   & 3.388 \\ \hline
7       & 1000 & 500 & 0.5  & 4               & 4.292        & 1.693 & 0.5     & 5.435                         & 1.168                & 1   & 3.454 & 1.105                               & 2.478                      & 5.549  & 6.09   & 3.388 \\ \hline
8       & 1000 & 500 & 0.5  & 5               & 5.283        & 1.694 & 0.5     & 7.211                         & 1.167                & 1   & 3.454 & 1.079                               & 2.479                      & 5.549  & 6.09   & 3.388 \\ \hline
9       & 600  & 600 & 1    & 2               & 2.711        & 2.24  & 0.5     & 2.73                          & 1.979                & 2   & 6.397 & 3.112                               & 5.571                      & 7.973  & 9.245  & 3.541 \\ \hline
10      & 600  & 600 & 1    & 3               & 3.621        & 2.24  & 0.5     & 4.262                         & 1.979                & 2   & 6.397 & 2.57                                & 5.571                      & 7.973  & 9.245  & 3.541 \\ \hline
11      & 600  & 600 & 1    & 4               & 4.585        & 2.239 & 0.5     & 5.955                         & 1.979                & 2   & 6.397 & 2.38                                & 5.571                      & 7.973  & 9.245  & 3.541 \\ \hline
12      & 600  & 600 & 1    & 5               & 5.566        & 2.24  & 0.5     & 7.716                         & 1.979                & 2   & 6.397 & 2.285                               & 5.571                      & 7.973  & 9.245  & 3.541 \\ \hline
13      & 400  & 600 & 1.5  & 2               & 3.066        & 2.722 & 0.5     & 3.319                         & 2.749                & 3   & 8.987 & 5.419                               & 9.017                      & 10.084 & 11.272 & 5.791 \\ \hline
14      & 400  & 600 & 1.5  & 3               & 3.931        & 2.72  & 0.5     & 4.803                         & 2.75                 & 3   & 8.987 & 4.238                               & 9.015                      & 10.084 & 11.272 & 5.791 \\ \hline
15      & 400  & 600 & 1.5  & 4               & 4.877        & 2.721 & 0.5     & 6.479                         & 2.749                & 3   & 8.987 & 3.826                               & 9.017                      & 10.084 & 11.272 & 5.791 \\ \hline
16      & 400  & 600 & 1.5  & 5               & 5.848        & 2.722 & 0.5     & 8.238                         & 2.749                & 3   & 8.987 & 3.619                               & 9.016                      & 10.084 & 11.272 & 5.791 \\ \hline
\end{tabular}
}
\centering
\caption{
$H_1$: the values associated with the gap conditions under each simulation settting.
For a visual representation of this table, please refer also to Figure \ref{fig:graph_of_gap_conditions_H1}.
}
\label{table:gap_conditions_H1}
\end{table}

\begin{table}[h!]
\scalebox{0.70}{
\begin{tabular}{|c|c|c|c|c|c|c|c|c|c|c|c|c|c|c|c|c|}
\hline
   & $n$  & $p$ & $c$  & $\lambda_{r_0}$ & $\psi_{r_0}$ & $b$   & $\mu_H$ & $g(\frac{\psi_{r_0}}{\mu_H})$ & $g(\frac{b}{\mu_H})$ & AIC & BIC   & $2\kappa(\frac{\psi_{r_0}}{\mu_H})$ & $2\kappa(\frac{b}{\mu_H})$ & $\pca\ica$  & $\pcb\icb$  & $\pcc\icc$ \\ \hline
1  & 1000 & 250 & 0.25 & 2               & 2.17         & 1.226 & 0.5     & 1.871                         & 0.554                & 0.5 & 1.727 & 0.563                               & 1.006                      & 3.732  & 3.97   & 2.813 \\ \hline
2  & 1000 & 250 & 0.25 & 3               & 3.151        & 1.226 & 0.5     & 3.463                         & 0.555                & 0.5 & 1.727 & 0.533                               & 1.005                      & 3.732  & 3.97   & 2.813 \\ \hline
3  & 1000 & 250 & 0.25 & 4               & 4.144        & 1.226 & 0.5     & 5.173                         & 0.554                & 0.5 & 1.727 & 0.522                               & 1.01                       & 3.732  & 3.97   & 2.813 \\ \hline
4  & 1000 & 250 & 0.25 & 5               & 5.14         & 1.227 & 0.5     & 6.949                         & 0.556                & 0.5 & 1.727 & 0.517                               & 1.004                      & 3.732  & 3.97   & 2.813 \\ \hline
5  & 1000 & 500 & 0.5  & 2               & 2.339        & 1.54  & 0.5     & 2.136                         & 0.955                & 1   & 3.454 & 1.24                                & 2.479                      & 5.549  & 6.09   & 3.388 \\ \hline
6  & 1000 & 500 & 0.5  & 3               & 3.303        & 1.54  & 0.5     & 3.718                         & 0.955                & 1   & 3.454 & 1.127                               & 2.483                      & 5.549  & 6.09   & 3.388 \\ \hline
7  & 1000 & 500 & 0.5  & 4               & 4.288        & 1.54  & 0.5     & 5.425                         & 0.955                & 1   & 3.454 & 1.086                               & 2.483                      & 5.549  & 6.09   & 3.388 \\ \hline
8  & 1000 & 500 & 0.5  & 5               & 5.279        & 1.54  & 0.5     & 7.203                         & 0.955                & 1   & 3.454 & 1.064                               & 2.482                      & 5.549  & 6.09   & 3.388 \\ \hline
9  & 600  & 600 & 1    & 2               & 2.678        & 2.077 & 0.5     & 2.678                         & 1.73                 & 2   & 6.397 & 2.939                               & 5.719                      & 7.973  & 9.245  & 3.541 \\ \hline
10 & 600  & 600 & 1    & 3               & 3.606        & 2.077 & 0.5     & 4.236                         & 1.73                 & 2   & 6.397 & 2.498                               & 5.718                      & 7.973  & 9.245  & 3.541 \\ \hline
11 & 600  & 600 & 1    & 4               & 4.575        & 2.077 & 0.5     & 5.936                         & 1.73                 & 2   & 6.397 & 2.336                               & 5.719                      & 7.973  & 9.245  & 3.541 \\ \hline
12 & 600  & 600 & 1    & 5               & 5.558        & 2.077 & 0.5     & 7.707                         & 1.73                 & 2   & 6.397 & 2.253                               & 5.72                       & 7.973  & 9.245  & 3.541 \\ \hline
13 & 400  & 600 & 1.5  & 2               & 3.018        & 2.552 & 0.5     & 3.237                         & 2.474                & 3   & 8.987 & 5.095                               & 9.301                      & 10.084 & 11.272 & 5.791 \\ \hline
14 & 400  & 600 & 1.5  & 3               & 3.908        & 2.552 & 0.5     & 4.761                         & 2.474                & 3   & 8.987 & 4.113                               & 9.301                      & 10.084 & 11.272 & 5.791 \\ \hline
15 & 400  & 600 & 1.5  & 4               & 4.863        & 2.552 & 0.5     & 6.45                          & 2.474                & 3   & 8.987 & 3.752                               & 9.303                      & 10.084 & 11.272 & 5.791 \\ \hline
16 & 400  & 600 & 1.5  & 5               & 5.837        & 2.552 & 0.5     & 8.219                         & 2.474                & 3   & 8.987 & 3.566                               & 9.304                      & 10.084 & 11.272 & 5.791 \\ \hline
\end{tabular}
}
\centering
\caption{
$H_2$: the values associated with the gap conditions under each simulation settting.
For a visual representation of this table, please refer also to Figure \ref{fig:graph_of_gap_conditions_H2}.
}
\label{table:gap_conditions_H2}
\end{table}

\begin{table}[h!]
\scalebox{0.70}{
\begin{tabular}{|c|c|c|c|c|c|c|c|c|c|c|c|c|c|c|c|c|}
\hline
   & $n$  & $p$ & $c$  & $\lambda_{r_0}$ & $\psi_{r_0}$ & $b$   & $\mu_H$ & $g(\frac{\psi_{r_0}}{\mu_H})$ & $g(\frac{b}{\mu_H})$ & AIC & BIC   & $2\kappa(\frac{\psi_{r_0}}{\mu_H})$ & $2\kappa(\frac{b}{\mu_H})$ & $\pca\ica$  & $\pcb\icb$  & $\pcc\icc$ \\ \hline
1  & 1000 & 250 & 0.25 & 2               & 2.209        & 1.666 & 0.501   & 1.926                         & 1.123                & 0.5 & 1.727 & 0.707                               & 1.502                      & 3.732  & 3.97   & 2.813 \\ \hline
2  & 1000 & 250 & 0.25 & 3               & 3.169        & 1.666 & 0.502   & 3.476                         & 1.123                & 0.5 & 1.727 & 0.598                               & 1.502                      & 3.732  & 3.97   & 2.813 \\ \hline
3  & 1000 & 250 & 0.25 & 4               & 4.155        & 1.665 & 0.501   & 5.18                          & 1.124                & 0.5 & 1.727 & 0.564                               & 1.502                      & 3.732  & 3.97   & 2.813 \\ \hline
4  & 1000 & 250 & 0.25 & 5               & 5.148        & 1.666 & 0.501   & 6.945                         & 1.124                & 0.5 & 1.727 & 0.547                               & 1.501                      & 3.732  & 3.97   & 2.813 \\ \hline
5  & 1000 & 500 & 0.5  & 2               & 2.419        & 2.033 & 0.501   & 2.252                         & 1.655                & 1   & 3.454 & 1.589                               & 2.925                      & 5.549  & 6.09   & 3.388 \\ \hline
6  & 1000 & 500 & 0.5  & 3               & 3.339        & 2.034 & 0.501   & 3.765                         & 1.657                & 1   & 3.454 & 1.272                               & 2.926                      & 5.549  & 6.09   & 3.388 \\ \hline
7  & 1000 & 500 & 0.5  & 4               & 4.31         & 2.033 & 0.501   & 5.444                         & 1.655                & 1   & 3.454 & 1.175                               & 2.925                      & 5.549  & 6.09   & 3.388 \\ \hline
8  & 1000 & 500 & 0.5  & 5               & 5.296        & 2.032 & 0.501   & 7.213                         & 1.655                & 1   & 3.454 & 1.129                               & 2.925                      & 5.549  & 6.09   & 3.388 \\ \hline
9  & 600  & 600 & 1    & 2               & 2.837        & 2.629 & 0.501   & 2.928                         & 2.59                 & 2   & 6.397 & 3.87                                & 5.993                      & 7.973  & 9.245  & 3.541 \\ \hline
10 & 600  & 600 & 1    & 3               & 3.677        & 2.628 & 0.501   & 4.348                         & 2.591                & 2   & 6.397 & 2.846                               & 5.993                      & 7.973  & 9.245  & 3.541 \\ \hline
11 & 600  & 600 & 1    & 4               & 4.621        & 2.631 & 0.502   & 5.993                         & 2.589                & 2   & 6.397 & 2.542                               & 5.995                      & 7.973  & 9.245  & 3.541 \\ \hline
12 & 600  & 600 & 1    & 5               & 5.592        & 2.63  & 0.501   & 7.743                         & 2.589                & 2   & 6.397 & 2.397                               & 5.994                      & 7.973  & 9.245  & 3.541 \\ \hline
13 & 400  & 600 & 1.5  & 2               & 3.256        & 3.15  & 0.501   & 3.624                         & 3.444                & 3   & 8.987 & 6.849                               & 9.329                      & 10.084 & 11.272 & 5.791 \\ \hline
14 & 400  & 600 & 1.5  & 3               & 4.017        & 3.15  & 0.501   & 4.929                         & 3.444                & 3   & 8.987 & 4.726                               & 9.33                       & 10.084 & 11.272 & 5.791 \\ \hline
15 & 400  & 600 & 1.5  & 4               & 4.931        & 3.149 & 0.501   & 6.553                         & 3.444                & 3   & 8.987 & 4.1                                 & 9.327                      & 10.084 & 11.272 & 5.791 \\ \hline
16 & 400  & 600 & 1.5  & 5               & 5.888        & 3.151 & 0.501   & 8.28                          & 3.445                & 3   & 8.987 & 3.805                               & 9.33                       & 10.084 & 11.272 & 5.791 \\ \hline
\end{tabular}
}
\centering
\caption{
$H_3$: the values associated with the gap conditions under each simulation settting.
For a visual representation of this table, please refer also to Figure \ref{fig:graph_of_gap_conditions_H3}.
}
\label{table:gap_conditions_H3}
\end{table}

%
%
\subsubsection{Results}\label{subsubsec:study1_results}

Tables \ref{table:unified_theorem_simulation_H1} to \ref{table:unified_theorem_simulation_H3} provide a detailed evaluation of each estimator's satisfaction of the gap conditions and their estimation accuracy across various simulation settings.
In these tables, ``T'' represents cases where the method satisfies the gap conditions, whereas ``F'' indicates failure to do so.
The values in the AIC to $\icc$ columns denote the percentage of correctly estimated ranks, based on $T=500$ repeated runs.

The estimators $\pca, \pcb$, and $\pcc$ require the estimation of $\sigma^2$ as mentioned in Section~\ref{subsec:intro_pc_ic}.
\citet{bai_ng2002} propose to estimate it by $\frac{1}{p-q} \sum_{j=q+1}^p \widehat{\lambda}_j$, however this is dependent on the upper bound of the search range for $r$, denoted as $q$.
In our preliminary experiments, we found that the choice of $q$ could lead to an underestimation of $\sigma^2$.
Therefore, in this simulation, we directly substitute the value of $\mu_H$ for $\widehat{\sigma}_q^2$ used in $\pca, \pcb$, and $\pcc$.

Tables \ref{table:unified_theorem_simulation_H1} to \ref{table:unified_theorem_simulation_H3} demonstrate that the estimation accuracy is generally high when the gap conditions are satisfied and low when they are not.
This observation substantiates our claim that under the general RMT assumptions (C1)-(C3), the selection consistency of information criterion-based estimators is entirely governed by the gap conditions.
However, it is worth noting that there are situations where the results do not align perfectly with the theory, especially when the gap conditions are marginally satisfied.

For instance, consider setting 15 under $H_2$ as detailed in Table~\ref{table:unified_theorem_simulation_H2}.
While both $\pcc$ and $\icc$ satisfy the gap conditions for this setting, their accuracy is less than ideal: $\pcc$ registers a value of 0.426 and $\icc$ a value of 0.578.
As illustrated in the fourth graph of Figure~\ref{fig:graph_of_gap_conditions_H2}, which corresponds to settings 13-16 under $H_2$, the intersection between the solid red curve and the plots for $\pcc/\icc$ occurs when $\lambda_{r_0}$ exceeds 3.5.
Moreover, the red solid curve closely approaches the plots for $\pcc/\icc$ at $\lambda_{r_0}=4$, indicating a narrow margin for satisfying the first gap condition.
It is worth noting that the first gap condition ensures that the selection criterion reaches its minimum value at $r=r_0$ for $1 \leq r \leq r_0$ (see Appendix \ref{sec:proof_of_theorem} for a detailed proof).
This narrow margin suggests a potential risk of underestimating the rank.
In the context of setting 15 under $H_2$, the average value for $\widehat{r}_{\rm PC3}$ is 4.382, while for $\widehat{r}_{\rm IC3}$ it is 4.548, further reinforcing the tendency for slight underestimation in this specific scenario.

On the other hand, when viewed from a different angle, the following observations can be made:
In contrast to more straightforward tasks like binary classification, which afford a 50\% accuracy through random guessing, our current task requires the estimation of $r_0=5$ from among $q$ candidates.
Given the inherent complexity of this problem, achieving accuracy rates of 0.426 and 0.568 could, in a certain sense, be considered commendably high.
This level of accuracy may serve as empirical support for the validity of the unified selection consistency theorem.
Taking into account that computer simulations are constrained by finite $(n, p)$ values, and that generated data will inherently contain stochastic variations, expecting a sharp increase in accuracy precisely when $\lambda_{r_0}$ satisfies the first gap condition might be unrealistic.
To overcome the inherent randomness in the data and achieve a satisfactory level of accuracy in estimation, it appears that $\lambda_{r_0}$ needs to substantially surpass the threshold.

Similarly, GIC's estimation accuracy remains low in settings 1-4 under $H_3$ in Table~\ref{table:unified_theorem_simulation_H3}, despite satisfying the gap conditions.
As shown in the first graph of Figure~\ref{fig:graph_of_gap_conditions_H3}, corresponding to settings 1-4 under $H_3$, the dashed blue line, representing the second gap condition value, only slightly surpasses and nearly coincides with the solid blue line, representing the maximum noise strength.
Similar to the previous case, in this instance as well, the second gap condition was not met with a sufficient margin, and consequently, it was unable to deal with the randomness in the generated data, failing to reach a satisfactorily high level of accurate estimation.
The failure to satisfy the second gap condition may lead to an overestimation of the rank, in contrast to the violation of the first gap condition.
This is further supported by examining the computed $\widehat{r}_{\rm GIC}$ values in settings 1-4 under $H_3$ in Table~\ref{table:unified_theorem_simulation_H3}, which have an average around 19 to 21, substantially higher than the true rank $r_0=5$.

\begin{table}[h!]
\scalebox{0.70}{
\begin{tabular}{|c|c|c|c|c|c|c|c|c|c|c|c|c|c|}
\hline
   & $n$  & $p$ & $c$  & $\lambda_{r_0}$ & AIC       & BIC       & GIC       & $\pca$    & $\ica$    & $\pcb$    & $\icb$    & $\pcc$    & $\icc$    \\ \hline
1  & 1000 & 250 & 0.25 & 2               & 0.000 (F) & 0.716 (T) & 1.000 (T) & 0.000 (F) & 0.000 (F) & 0.000 (F) & 0.000 (F) & 0.000 (F) & 0.000 (F) \\ \hline
2  & 1000 & 250 & 0.25 & 3               & 0.000 (F) & 1.000 (T) & 1.000 (T) & 0.008 (F) & 0.042 (F) & 0.000 (F) & 0.002 (F) & 0.992 (T) & 0.994 (T) \\ \hline
3  & 1000 & 250 & 0.25 & 4               & 0.000 (F) & 1.000 (T) & 0.996 (T) & 1.000 (T) & 1.000 (T) & 1.000 (T) & 1.000 (T) & 1.000 (T) & 1.000 (T) \\ \hline
4  & 1000 & 250 & 0.25 & 5               & 0.000 (F) & 1.000 (T) & 0.998 (T) & 1.000 (T) & 1.000 (T) & 1.000 (T) & 1.000 (T) & 1.000 (T) & 1.000 (T) \\ \hline
5  & 1000 & 500 & 0.5  & 2               & 0.000 (F) & 0.000 (F) & 1.000 (T) & 0.000 (F) & 0.000 (F) & 0.000 (F) & 0.000 (F) & 0.000 (F) & 0.000 (F) \\ \hline
6  & 1000 & 500 & 0.5  & 3               & 0.000 (F) & 0.676 (T) & 1.000 (T) & 0.000 (F) & 0.000 (F) & 0.000 (F) & 0.000 (F) & 0.728 (T) & 0.804 (T) \\ \hline
7  & 1000 & 500 & 0.5  & 4               & 0.000 (F) & 1.000 (T) & 1.000 (T) & 0.042 (F) & 0.086 (F) & 0.000 (F) & 0.000 (F) & 1.000 (T) & 1.000 (T) \\ \hline
8  & 1000 & 500 & 0.5  & 5               & 0.000 (F) & 1.000 (T) & 1.000 (T) & 1.000 (T) & 1.000 (T) & 0.980 (T) & 0.980 (T) & 1.000 (T) & 1.000 (T) \\ \hline
9  & 600  & 600 & 1    & 2               & 0.838 (T) & 0.000 (F) & 0.124 (F) & 0.000 (F) & 0.000 (F) & 0.000 (F) & 0.000 (F) & 0.000 (F) & 0.000 (F) \\ \hline
10 & 600  & 600 & 1    & 3               & 0.836 (T) & 0.000 (F) & 1.000 (T) & 0.000 (F) & 0.000 (F) & 0.000 (F) & 0.000 (F) & 0.958 (T) & 0.978 (T) \\ \hline
11 & 600  & 600 & 1    & 4               & 0.806 (T) & 0.002 (F) & 1.000 (T) & 0.000 (F) & 0.000 (F) & 0.000 (F) & 0.000 (F) & 1.000 (T) & 1.000 (T) \\ \hline
12 & 600  & 600 & 1    & 5               & 0.808 (T) & 0.954 (T) & 1.000 (T) & 0.006 (F) & 0.026 (F) & 0.000 (F) & 0.000 (F) & 1.000 (T) & 1.000 (T) \\ \hline
13 & 400  & 600 & 1.5  & 2               & 0.828 (T) & 0.000 (F) & 0.000 (F) & 0.000 (F) & 0.000 (F) & 0.000 (F) & 0.000 (F) & 0.000 (F) & 0.000 (F) \\ \hline
14 & 400  & 600 & 1.5  & 3               & 0.998 (T) & 0.000 (F) & 0.868 (T) & 0.000 (F) & 0.000 (F) & 0.000 (F) & 0.000 (F) & 0.000 (F) & 0.000 (F) \\ \hline
15 & 400  & 600 & 1.5  & 4               & 1.000 (T) & 0.000 (F) & 1.000 (T) & 0.000 (F) & 0.000 (F) & 0.000 (F) & 0.000 (F) & 0.442 (T) & 0.588 (T) \\ \hline
16 & 400  & 600 & 1.5  & 5               & 1.000 (T) & 0.000 (F) & 1.000 (T) & 0.000 (F) & 0.000 (F) & 0.000 (F) & 0.000 (F) & 1.000 (T) & 1.000 (T) \\ \hline
\end{tabular}
}
\centering
\caption{
$H_1$: The estimation accuracy rates for each estimator in each setting, based on 500 repeated runs.
The values within each cell represent the computed estimation accuracy rates, and the adjacent T or F indicates whether the gap conditions for that estimator are satisfied in the given setting.
}
\label{table:unified_theorem_simulation_H1}
\end{table}

\begin{table}[h!]
\scalebox{0.70}{
\begin{tabular}{|c|c|c|c|c|c|c|c|c|c|c|c|c|c|}
\hline
   & $n$  & $p$ & $c$  & $\lambda_{r_0}$ & AIC       & BIC       & GIC       & $\pca$    & $\ica$    & $\pcb$    & $\icb$    & $\pcc$    & $\icc$    \\ \hline
1  & 1000 & 250 & 0.25 & 2               & 0.134 (F) & 0.650 (T) & 1.000 (T) & 0.000 (F) & 0.000 (F) & 0.000 (F) & 0.000 (F) & 0.000 (F) & 0.000 (F) \\ \hline
2  & 1000 & 250 & 0.25 & 3               & 0.094 (F) & 1.000 (T) & 1.000 (T) & 0.002 (F) & 0.020 (F) & 0.000 (F) & 0.000 (F) & 0.984 (T) & 0.994 (T) \\ \hline
3  & 1000 & 250 & 0.25 & 4               & 0.110 (F) & 1.000 (T) & 1.000 (T) & 1.000 (T) & 1.000 (T) & 1.000 (T) & 1.000 (T) & 1.000 (T) & 1.000 (T) \\ \hline
4  & 1000 & 250 & 0.25 & 5               & 0.130 (F) & 1.000 (T) & 1.000 (T) & 1.000 (T) & 1.000 (T) & 1.000 (T) & 1.000 (T) & 1.000 (T) & 1.000 (T) \\ \hline
5  & 1000 & 500 & 0.5  & 2               & 0.984 (T) & 0.000 (F) & 1.000 (T) & 0.000 (F) & 0.000 (F) & 0.000 (F) & 0.000 (F) & 0.000 (F) & 0.000 (F) \\ \hline
6  & 1000 & 500 & 0.5  & 3               & 0.990 (T) & 0.638 (T) & 1.000 (T) & 0.000 (F) & 0.000 (F) & 0.000 (F) & 0.000 (F) & 0.690 (T) & 0.774 (T) \\ \hline
7  & 1000 & 500 & 0.5  & 4               & 0.986 (T) & 1.000 (T) & 1.000 (T) & 0.028 (F) & 0.056 (F) & 0.000 (F) & 0.000 (F) & 1.000 (T) & 1.000 (T) \\ \hline
8  & 1000 & 500 & 0.5  & 5               & 0.978 (T) & 1.000 (T) & 1.000 (T) & 1.000 (T) & 1.000 (T) & 0.964 (T) & 0.974 (T) & 1.000 (T) & 1.000 (T) \\ \hline
9  & 600  & 600 & 1    & 2               & 1.000 (T) & 0.000 (F) & 0.298 (F) & 0.000 (F) & 0.000 (F) & 0.000 (F) & 0.000 (F) & 0.000 (F) & 0.000 (F) \\ \hline
10 & 600  & 600 & 1    & 3               & 1.000 (T) & 0.000 (F) & 1.000 (T) & 0.000 (F) & 0.000 (F) & 0.000 (F) & 0.000 (F) & 0.906 (T) & 0.952 (T) \\ \hline
11 & 600  & 600 & 1    & 4               & 1.000 (T) & 0.004 (F) & 1.000 (T) & 0.000 (F) & 0.000 (F) & 0.000 (F) & 0.000 (F) & 1.000 (T) & 1.000 (T) \\ \hline
12 & 600  & 600 & 1    & 5               & 1.000 (T) & 0.944 (T) & 1.000 (T) & 0.004 (F) & 0.006 (F) & 0.000 (F) & 0.000 (F) & 1.000 (T) & 1.000 (T) \\ \hline
13 & 400  & 600 & 1.5  & 2               & 0.640 (T) & 0.000 (F) & 0.000 (F) & 0.000 (F) & 0.000 (F) & 0.000 (F) & 0.000 (F) & 0.000 (F) & 0.000 (F) \\ \hline
14 & 400  & 600 & 1.5  & 3               & 1.000 (T) & 0.000 (F) & 0.922 (T) & 0.000 (F) & 0.000 (F) & 0.000 (F) & 0.000 (F) & 0.000 (F) & 0.000 (F) \\ \hline
15 & 400  & 600 & 1.5  & 4               & 1.000 (T) & 0.000 (F) & 1.000 (T) & 0.000 (F) & 0.000 (F) & 0.000 (F) & 0.000 (F) & 0.426 (T) & 0.578 (T) \\ \hline
16 & 400  & 600 & 1.5  & 5               & 1.000 (T) & 0.002 (F) & 1.000 (T) & 0.000 (F) & 0.000 (F) & 0.000 (F) & 0.000 (F) & 0.998 (T) & 0.998 (T) \\ \hline
\end{tabular}
}
\centering
\caption{
$H_2$: The estimation accuracy rates for each estimator in each setting, based on 500 repeated runs.
The values within each cell represent the computed estimation accuracy rates, and the adjacent T or F indicates whether the gap conditions for that estimator are satisfied in the given setting.
}
\label{table:unified_theorem_simulation_H2}
\end{table}

\begin{table} [hbt!]
\scalebox{0.70}{
\begin{tabular}{|c|c|c|c|c|c|c|c|c|c|c|c|c|c|}
\hline
   & $n$  & $p$ & $c$  & $\lambda_{r_0}$ & AIC       & BIC       & GIC       & $\pca$    & $\ica$    & $\pcb$    & $\icb$    & $\pcc$    & $\icc$    \\ \hline
1  & 1000 & 250 & 0.25 & 2               & 0.000 (F) & 0.770 (T) & 0.326 (T) & 0.000 (F) & 0.000 (F) & 0.000 (F) & 0.000 (F) & 0.000 (F) & 0.000 (F) \\ \hline
2  & 1000 & 250 & 0.25 & 3               & 0.000 (F) & 1.000 (T) & 0.384 (T) & 0.000 (F) & 0.084 (F) & 0.000 (F) & 0.024 (F) & 0.988 (T) & 0.972 (T) \\ \hline
3  & 1000 & 250 & 0.25 & 4               & 0.000 (F) & 1.000 (T) & 0.404 (T) & 1.000 (T) & 1.000 (T) & 1.000 (T) & 0.994 (T) & 1.000 (T) & 1.000 (T) \\ \hline
4  & 1000 & 250 & 0.25 & 5               & 0.000 (F) & 1.000 (T) & 0.372 (T) & 1.000 (T) & 1.000 (T) & 1.000 (T) & 1.000 (T) & 1.000 (T) & 1.000 (T) \\ \hline
5  & 1000 & 500 & 0.5  & 2               & 0.000 (F) & 0.000 (F) & 0.994 (T) & 0.000 (F) & 0.000 (F) & 0.000 (F) & 0.000 (F) & 0.000 (F) & 0.000 (F) \\ \hline
6  & 1000 & 500 & 0.5  & 3               & 0.000 (F) & 0.716 (T) & 1.000 (T) & 0.000 (F) & 0.000 (F) & 0.000 (F) & 0.000 (F) & 0.758 (T) & 0.822 (T) \\ \hline
7  & 1000 & 500 & 0.5  & 4               & 0.000 (F) & 1.000 (T) & 1.000 (T) & 0.036 (F) & 0.134 (F) & 0.000 (F) & 0.004 (F) & 1.000 (T) & 1.000 (T) \\ \hline
8  & 1000 & 500 & 0.5  & 5               & 0.000 (F) & 1.000 (T) & 1.000 (T) & 0.998 (T) & 1.000 (T) & 0.962 (T) & 0.958 (T) & 1.000 (T) & 1.000 (T) \\ \hline
9  & 600  & 600 & 1    & 2               & 0.000 (F) & 0.000 (F) & 0.016 (F) & 0.000 (F) & 0.000 (F) & 0.000 (F) & 0.000 (F) & 0.000 (F) & 0.000 (F) \\ \hline
10 & 600  & 600 & 1    & 3               & 0.000 (F) & 0.000 (F) & 1.000 (T) & 0.000 (F) & 0.000 (F) & 0.000 (F) & 0.000 (F) & 0.976 (T) & 0.984 (T) \\ \hline
11 & 600  & 600 & 1    & 4               & 0.000 (F) & 0.014 (F) & 1.000 (T) & 0.000 (F) & 0.000 (F) & 0.000 (F) & 0.000 (F) & 1.000 (T) & 1.000 (T) \\ \hline
12 & 600  & 600 & 1    & 5               & 0.000 (F) & 0.942 (T) & 1.000 (T) & 0.000 (F) & 0.028 (F) & 0.000 (F) & 0.000 (F) & 1.000 (T) & 1.000 (T) \\ \hline
13 & 400  & 600 & 1.5  & 2               & 0.000 (F) & 0.000 (F) & 0.000 (F) & 0.000 (F) & 0.000 (F) & 0.000 (F) & 0.000 (F) & 0.000 (F) & 0.000 (F) \\ \hline
14 & 400  & 600 & 1.5  & 3               & 0.000 (F) & 0.000 (F) & 0.644 (T) & 0.000 (F) & 0.000 (F) & 0.000 (F) & 0.000 (F) & 0.000 (F) & 0.000 (F) \\ \hline
15 & 400  & 600 & 1.5  & 4               & 0.000 (F) & 0.000 (F) & 1.000 (T) & 0.000 (F) & 0.000 (F) & 0.000 (F) & 0.000 (F) & 0.474 (T) & 0.680 (T) \\ \hline
16 & 400  & 600 & 1.5  & 5               & 0.000 (F) & 0.002 (F) & 1.000 (T) & 0.000 (F) & 0.000 (F) & 0.000 (F) & 0.000 (F) & 1.000 (T) & 1.000 (T) \\ \hline
\end{tabular}
}
\centering
\caption{
$H_3$: The estimation accuracy rates for each estimator in each setting, based on 500 repeated runs.
The values within each cell represent the computed estimation accuracy rates, and the adjacent T or F indicates whether the gap conditions for that estimator are satisfied in the given setting.
}
\label{table:unified_theorem_simulation_H3}
\end{table}

%
%

\subsection{Study II: exploration of characteristics in information criterion-based estimators}\label{subsec:study2}

In this section, we aim to delve deeper into the properties of information criterion-based estimators, such as understanding their advantages and limitations.
Furthermore, we have also selected some methods from Categories 2 and 3 to compare the performance of information criterion-based estimators from different perspectives.
However, please be informed that since the focus of this research is on information criterion-based methods, we will not delve deeply into the characteristics of the methods from Categories 2 and 3.
When devising the simulation settings for {\it Study II}, we paid attention to the fact that as long as the general RMT assumptions (C1)-(C3) are satisfied, the following conditions would not influence selection consistency.

\begin{itemize}\itemsep=0pt
\item A high rank $r_0$ is set.
\item The entries of $\bm{z}_i \eqdef (z_{i,1}, \ldots, z_{i,p})^\top$ are generated from a heavy-tailed distribution.
\item The population covariance matrix $\bm{\Sigma}$ is sparse.
\item The largest population eigenvalue $\lambda_1$ is significantly larger than $\lambda_2$.
\item The signal eigenvalues $\lambda_{r_0}, \ldots, \lambda_{1}$ form a geometric progression with a common ratio greater than 1.
\end{itemize}

In Sections \ref{subsubsec:study2_high_rank}-\ref{subsubsec:study2_power_spiked}, we closely adhere to the simulation procedure (Steps 1 to 5) outlined in Section~\ref{subsubsec:study1_settings}.
However, we make slight adjustments:
\begin{itemize}\itemsep=0pt
\item In Section \ref{subsubsec:study2_high_rank}, we vary the value of $r_0$ in each individual simulation, which was originally fixed at $5$.
(In Sections \ref{subsubsec:study2_heavy_tailed} through \ref{subsubsec:study2_power_spiked}, we reset $r_0$ to 5.)
\item In Section \ref{subsubsec:study2_heavy_tailed}, we replace the distribution of $\bm{x}_1, \ldots, \bm{x}_n \iid N(\bm{0},\bm{\Sigma})$ with other distributions.
\item In Section \ref{subsubsec:study2_sparse_covariance}, we modify the generation method of $\bm{\Gamma}$ to induce sparsity in $\bm{\Sigma}$.
\item In Section \ref{subsubsec:study2_large_lambda1}, we explicitly define the value of $\lambda_1$. (The method for generating $\lambda_2,\ldots,\lambda_p$ remains unchanged.)
\item In Section \ref{subsubsec:study2_power_spiked}, we alter the generation method of $\bm{\Lambda}$ to establish a geometric progression for $\lambda_{r_0}, \ldots, \lambda_{1}$.
\end{itemize}
The following condition remains consistent across Sections~\ref{subsubsec:study2_high_rank}-\ref{subsubsec:study2_power_spiked}:
\begin{itemize}\itemsep=0pt
  \item The sample size and the data dimension $(n,p)$ are set to $(500, 200)$ in all settings.
  \item The signal population eigenvalue, $\lambda_{r_0}$, is set to $6$, a value sufficiently large to ensure that $r_0$ is identified as the rank.
  \item The LSD of the population covariance matrix is set as $H(t) \eqdef \mathcal{I}(1 \leq t)$.
\end{itemize}

For comparison, we select four additional methods within Categories 2 and 3: ACT \cite{act}, DPA \citep{dpa}, ED \citep{ed} and GR \citep{er_gr}.
Though there are several other methods from Categories 2 and 3, we are particularly concerned that including results from all methods would make the tables overly complex.
Furthermore, the primary focus of this article is on the information criterion-based estimators.
Broadly evaluating the performance of all rank estimation methods would diverge from the main focus of this article.
Therefore, we have chosen two methods from Category 2 and two from Category 3 for comparison.
Through the simulations of the five different scenarios mentioned earlier, we attempt to investigate the advantages and limitations of information criterion-based estimators by comparing their performance with the selected alternative methods.
Due to space limitations, we also exclude the results of $\pca$, $\ica$, $\pcb$, and $\icb$.
As we also did in Section \ref{subsec:study1} ({\it Study I}), we do not estimate $\sigma^2$ in $\pcc$, opting instead to directly substitute $\widehat{\sigma}_q^2 \leftarrow \mu_H \stackrel{\rm eq}{=} 1$.

%
%
\subsubsection{Model with high rank $r_0$}\label{subsubsec:study2_high_rank}

The simulations in this study follow the procedure outlined in Section \ref{subsubsec:study1_settings}, with the value of $r_0$ varied to be $5$, $10$, $15$, $20$, and $25$.
The upper bound for the search range of variable $r$, denoted by $q$, is set to $\lfloor \sqrt{2\min\{n,p\}} \rfloor = 28$.
The accuracy of the estimates for each method in each setting is presented in Table \ref{table:study2_high_rank}.

Contrary to our expectations, we found that larger values of $r_0$ impact BIC, $\pcc$ and $\icc$.
Theoretically, the selection consistency of information criterion-based estimators should not be influenced by an increase in $r_0$.
However, in practice, finite values must be set for $(n,p)$ during numerical experiments, and results consistent with the theoretical conclusions under the assumption of $(n,p)$ diverging to infinity may not necessarily be observed.
In our specific setting with $p=200$, it appears that values such as $r_0=20$ and $r_0=25$ are too substantial to be ignored.
With the increase of $r_0$, there is a decrease in the smallest sample signal eigenvalue $\widehat{\lambda}_{r_0}$.
As a result, methods like BIC, $\pcc$ and $\icc$, which rely on larger signal eigenvalues, are more prone to this decrease.
Interestingly, despite sharing the same gap conditions, $\pcc$ appears to be more significantly impacted in terms of estimation accuracy in this scenario compared to $\icc$.
Recall that $\pcc$ and $\icc$ are defined by Equations (\ref{eq:pc3_estimator}) and (\ref{eq:ic3_estimator}), respectively, and we directly set $\widehat{\sigma}_q^2 \leftarrow \mu_H \stackrel{\rm eq}{=} 1$ in our simulations.
The primary difference between $\pcc$ and $\icc$ is the presence of natural logarithm in the first term, while the penalty terms are identical.
This means that $\pcc$ imposes a larger penalty than $\icc$, requiring stronger signal eigenvalues.
We attribute this to the reason why $\pcc$ is more significantly impacted than $\icc$.

\begin{table}[h!]
\scalebox{0.9}{
\begin{tabular}{|c|c|c|c|c|c|c|c|c|c|c|}
\hline
Setting & $r_0$ & AIC   & BIC   & GIC   & $\pcc$ & $\icc$ & ACT   & DPA   & ED    & GR    \\ \hline
1       & 5     & 1.000 & 1.000 & 1.000 & 0.988  & 0.998  & 1.000 & 1.000 & 0.986 & 1.000 \\ \hline
2       & 10    & 1.000 & 0.986 & 1.000 & 0.668  & 0.978  & 1.000 & 1.000 & 0.994 & 1.000 \\ \hline
3       & 15    & 1.000 & 0.916 & 1.000 & 0.068  & 0.910  & 1.000 & 1.000 & 0.990 & 1.000 \\ \hline
4       & 20    & 1.000 & 0.710 & 1.000 & 0.000  & 0.792  & 1.000 & 1.000 & 1.000 & 1.000 \\ \hline
5       & 25    & 0.998 & 0.444 & 1.000 & 0.000  & 0.686  & 1.000 & 0.236 & 1.000 & 1.000 \\ \hline
\end{tabular}
}
\centering
\caption{
Impact of large $r_0$.
Simulations follow the procedure in Section \ref{subsubsec:study1_settings}, with $r_0$ values of 10, 15, 20, and 25.
The upper bound for the seach range of $r$ is set to 28.
Larger $r_0$ results in decreased sample signal  eigenvalue $\widehat{\lambda}_{r_0}$, impacting the performance of BIC, $\pcc$ and $\icc$, which rely on a strong signal.
}
\label{table:study2_high_rank}
\end{table}

%
%
\subsubsection{Model with a heavy-tailed distribution}\label{subsubsec:study2_heavy_tailed}

In this section, we conduct simulations mostly following the procedure outlined in Section \ref{subsubsec:study1_settings}, with a modification to Step 3 in which we replace the distribution $N(\bm{0}, \bm{\Sigma})$ with alternative distributions.
We generate $\bm{x}_i \eqdef \bm{\Sigma}^{1/2} \bm{z}_i$, where the entries of $\bm{z}_i$ (denoted as $z_{i,j}$) independently follow heavy-tailed distributions with mean $0$ and variance $1$.
We generate $z_{i,j}$ using the formula $z_{i,j} \eqdef (y_{i,j} - {\rm E}[y_{i,j}])/ {\rm Var}[y_{i,j}]^{1/2}$, where the distribution of $y_{i,j}$ is as follows:
\begin{itemize}\itemsep=0pt
\item $t$-distribution with degrees of freedom $\nu=5$.
\item Pareto distribution with probability density function (p.d.f.) $f(y) \eqdef \frac{\alpha \nu^\alpha}{y^{\alpha+1}} \mathcal{I}(y \geq \nu)$, where $\alpha=5$ and $\nu=1$.
\item Log-normal distribution with $\mu=0$ and $\sigma^2=1$.
\end{itemize}

The accuracy of the estimates for each method in each setting is presented in Table \ref{table:study2_heavy_tailed}.
We observed that the accuracy of AIC and GIC decreased, while the accuracy for BIC, $\pcc$ and $\icc$ was less impacted.
We attribute this phenomenon to the following cause: when data are generated from a heavy-tailed distribution, the observed sample noise eigenvalue $\widehat{\lambda}_{r_0+1}$ tends to be larger than in the case of a normal distribution.
Therefore, methods that are robust to large noise (but require a strong signal), such as BIC, $\pcc$ and $\icc$, are less susceptible to this scenario, whereas methods that rely on small noise (but are sensitive to a small signal), such as AIC, are more affected.

While further simulations will continue in Sections \ref{subsubsec:study2_sparse_covariance}-\ref{subsubsec:study2_power_spiked}, we find it important to highlight some key observations at this juncture, taking into account the results of Sections \ref{subsubsec:study2_high_rank} and \ref{subsubsec:study2_heavy_tailed}.
In Section \ref{subsubsec:study2_high_rank}, we observed that as $r_0$ increased, the sample signal eigenvalue $\widehat{\lambda}_{r_0}$ decreased, adversely affecting methods that require strong signals like BIC, $\pcc$, and $\icc$.
In Section \ref{subsubsec:study2_heavy_tailed}, switching to a heavy-tailed distribution led to an increase in the sample noise eigenvalue $\widehat{\lambda}_{r_0+1}$, which adversely affected methods that require low noise levels, such as AIC.
For the information criterion-based estimators, whose the selection consistency is governed by the gap conditions, there exists a trade-off between sensitivity to signal and robustness to noise; they are likely to be affected in at least one of these scenarios detailed in Sections \ref{subsubsec:study2_high_rank} and \ref{subsubsec:study2_heavy_tailed}.
However, ACT consistently maintains high estimation accuracy in both scenarios, suggesting its potential to outperform information criterion-based estimators under specific conditions.
We note these phenomena as limitations of information criterion-based estimators.

\begin{table}[h!]
\scalebox{0.9}{
\begin{tabular}{|c|c|c|c|c|c|c|c|c|c|c|}
\hline
Setting & dist       & AIC   & BIC   & GIC   & $\pcc$ & $\icc$ & ACT   & DPA   & ED    & GR    \\ \hline
1       & $t_5$      & 0.914 & 1.000 & 0.976 & 0.984  & 1.000  & 0.998 & 0.970 & 0.894 & 0.994 \\ \hline
2       & Pareto     & 0.018 & 0.928 & 0.522 & 0.916  & 0.930  & 0.996 & 0.718 & 0.292 & 0.848 \\ \hline
3       & Log-normal & 0.000 & 0.860 & 0.206 & 0.882  & 0.882  & 0.988 & 0.698 & 0.164 & 0.756 \\ \hline
\end{tabular}
}
\centering
\caption{
Impact of switching to heavy-tailed distributions. ($r_0$ is fixed at 5.)
In this scenario, the entries of random vector $\bm{z}_i$ are independently and identically generated from heavy-tailed distributions, rather than a normal distribution.
The chosen distributions include a $t$-distribution with 5 degrees of freedom, a Pareto distribution, and a log-normal distribution.
These distributions are shifted and scaled to have a mean of 1 and variance of 0.
Shifting to a heavy-tailed distribution leads to an increase in the noise sample eigenvalue $\widehat{\lambda}_{r_0+1}$, which significantly affects the performance of AIC (and GIC) due to their reliance on low noise levels.
}
\label{table:study2_heavy_tailed}
\end{table}

%
%
\subsubsection{Model with a sparse population covariance matrix}\label{subsubsec:study2_sparse_covariance}

In this section, we modify Step 2 of the simulation procedure outlined in Section \ref{subsubsec:study1_settings} to induce sparsity in the population covariance matrix $\bm{\Sigma}$.
To achieve this, we generate $\bm{\Gamma}$ as a block diagonal matrix, ${\rm blockdiag}(\bm{\Gamma}_1, \ldots, \bm{\Gamma}_K)$, where $\bm{\Gamma}_k$ is a $p_k \times p_k$ orthogonal matrix drawn uniformly.
$K$ is the number of blocks, and we set the size of each block, denoted as $p_k$, to $\lfloor p/K \rfloor$ for $1 \leq k \leq K-1$, with $\sum_{\ell=1}^K p_\ell = p$.
Importantly, increasing the number of blocks results in a sparser $\bm{\Gamma}$, thereby inducing sparsity in $\bm{\Sigma}$ as well.
We set the number of blocks to $K=5, 10, 15$, and $20$.

The accuracy of the estimates for each method in each setting is presented in Table \ref{table:study2_sparse_covariance}.
None of the methods based on the information criterion were affected by the sparse covariance matrix.
However, methods that rely on the pervasiveness in the population covariance matrix, such as ACT, are adversely affected by the sparse covariance matrix.
\begin{table}[h!]
\scalebox{0.9}{
\begin{tabular}{|c|c|c|c|c|c|c|c|c|c|c|}
\hline
Setting & Block & AIC   & BIC   & GIC   & $\pcc$ & $\icc$ & ACT   & DPA   & ED    & GR    \\ \hline
1       & 5     & 1.000 & 1.000 & 1.000 & 0.984  & 0.996  & 0.982 & 1.000 & 0.976 & 1.000 \\ \hline
2       & 10    & 1.000 & 1.000 & 1.000 & 0.986  & 0.998  & 0.840 & 1.000 & 0.978 & 1.000 \\ \hline
3       & 15    & 1.000 & 1.000 & 1.000 & 0.986  & 0.998  & 0.012 & 0.918 & 0.974 & 1.000 \\ \hline
4       & 20    & 1.000 & 0.996 & 1.000 & 0.978  & 0.994  & 0.000 & 0.372 & 0.976 & 1.000 \\ \hline
\end{tabular}
}
\centering
\caption{
Impact of sparse population covariance matrix. ($r_0$ is fixed at $5$.)
The population covariance matrix, denoted as $\bm{\Sigma}$, is intentionally sparse.
Information criterion-based estmators remain unaffected by the sparseness of the population covariance matrix.
However, it is worth noting that certain other methods may be influenced by this change in sparsity.
}
\label{table:study2_sparse_covariance}
\end{table}

%
%
\subsubsection{Model with extremely large $\lambda_1$}\label{subsubsec:study2_large_lambda1}

In this section, we incorporate an extremely large population eigenvalue $\lambda_1 \gg \lambda_2$ into the generation method of $\bm{\Lambda}$, deviating from the original Step 1 outlined in Section \ref{subsubsec:study1_settings}.
Specifically, we set $\lambda_1$ to be $10^2$, $10^3$, $10^4$, $10^5$ and $10^6$, while generating $\lambda_2, \ldots, \lambda_{r_0-1}$ using $U_j + \lambda_{r_0}$ with $U_j \iid U(0,1)$, following the procedure outlined in Section \ref{subsubsec:study1_settings}.

The accuracy of the estimates for each method in each setting is presented in Table \ref{table:study2_large_lambda1}.
The information criterion-based estimators are not affected by the extremely large $\lambda_1$.
However, ACT, DPA, and GR show decreased estimation accuracy when confronted with extremely large $\lambda_1$.

\begin{table}[h!]
\scalebox{0.9}{
\begin{tabular}{|c|c|c|c|c|c|c|c|c|c|c|}
\hline
Setting & $\lambda_1$ & AIC   & BIC   & GIC   & $\pcc$ & $\icc$ & ACT   & DPA   & ED    & GR    \\ \hline
1       & $10^2$      & 1.000 & 1.000 & 1.000 & 0.990  & 1.000  & 1.000 & 0.428 & 0.972 & 0.000 \\ \hline
2       & $10^3$      & 1.000 & 1.000 & 1.000 & 0.992  & 0.994  & 1.000 & 0.000 & 0.972 & 0.000 \\ \hline
3       & $10^4$      & 1.000 & 1.000 & 1.000 & 0.992  & 0.998  & 0.000 & 0.000 & 0.988 & 0.000 \\ \hline
4       & $10^5$      & 1.000 & 0.998 & 1.000 & 0.986  & 0.992  & 0.000 & 0.000 & 0.978 & 0.000 \\ \hline
5       & $10^6$      & 1.000 & 1.000 & 1.000 & 0.996  & 1.000  & 0.000 & 0.000 & 0.974 & 0.000 \\ \hline
\end{tabular}
}
\centering
\caption{
Impact of extremely large $\lambda_1$. ($r_0$ is fixed at $5$.)
In this investigation, the largest population eigenvalue, $\lambda_1$, is explicitly set to five different values: $10^2$, $10^3$, $10^4$, $10^5$ and $10^6$.
Notably, information criterion-based estimators demonstrate consistent performance regardless of this change.
However, it is important to recognize that certain other methods may exhibit vulnerability to this alteration.
}
\label{table:study2_large_lambda1}
\end{table}

%
%
\subsubsection{Power spiked model}\label{subsubsec:study2_power_spiked}

In this section, we modify the original Step 1 as outlined in Section \ref{subsubsec:study1_settings} to alter the generation method of $\bm{\Lambda}$.
Our modification creates a power spiked model where $\lambda_1, \ldots, \lambda_{r_0}$ form a geometric progression.
We define the common ratio as $\alpha > 1$ and set $\lambda_1 = \lambda_{r_0} \alpha^{r_0-1}$, $\lambda_2 = \lambda_{r_0} \alpha^{r_0-2}$, $\ldots$, $\lambda_{r_0-1} = \lambda_{r_0} \alpha$.
We consider four different values for the common ratio: $\alpha = 2, 3, 4,$ and $5$.

The accuracy of the estimates for each method in each setting is presented in Table \ref{table:study2_power_spiked}.
The information criterion-based estimators are not influenced by the power spiked model.
However, other methods such as ACT, DPA, and GR exhibit sensitivity to the power spiked model.

\begin{table}[h!]
\scalebox{0.9}{
\begin{tabular}{|c|c|c|c|c|c|c|c|c|c|c|}
\hline
Setting & $\alpha$ & AIC   & BIC   & GIC   & $\pcc$ & $\icc$ & ACT   & DPA   & ED    & GR    \\ \hline
1       & 2        & 1.000 & 1.000 & 1.000 & 1.000  & 1.000  & 1.000 & 0.250 & 0.990 & 0.960 \\ \hline
2       & 3        & 1.000 & 1.000 & 1.000 & 1.000  & 1.000  & 0.980 & 0.000 & 0.980 & 0.150 \\ \hline
3       & 4        & 1.000 & 1.000 & 1.000 & 1.000  & 1.000  & 0.000 & 0.000 & 0.970 & 0.000 \\ \hline
4       & 5        & 1.000 & 1.000 & 1.000 & 1.000  & 1.000  & 0.000 & 0.000 & 0.970 & 0.000 \\ \hline
\end{tabular}
}
\centering
\caption{
Impact of power spiked model. ($r_0$ is fixed at $5$.)
In this analysis, the signal population eigenvalues $\lambda_1, \ldots, \lambda_{r_0}$ are intentionally arranged in a geometric progression, forming a power spiked model.
Notably, information criterion-based estimators demonstrate robustness in this model configuration.
However, it is important to note that certain other methods exhibit vulnerability to the power spiked model.
}
\label{table:study2_power_spiked}
\end{table}

\subsubsection{Summary}

In summarizing the simulation results presented in Sections \ref{subsubsec:study2_high_rank} to \ref{subsubsec:study2_power_spiked}, we can outline the characteristics of the information criterion-based estimators as follows:

\begin{itemize}\itemsep=0pt
\item
In Sections \ref{subsubsec:study2_high_rank} and \ref{subsubsec:study2_heavy_tailed}, we find some unexpected results in our numerical experiments.
Although theory suggests that the value of $r_0$ and the transition from a normal to a heavy-tailed distribution should not impact the consistency of estimators based on the information criterion, our findings contradict this expectation.
Specifically, an increase in $r_0$ relative to the given $p$ leads to a decrease in the sample signal eigenvalue $ \widehat{\lambda}_{r_0}$, thereby affecting estimators that rely on a strong signal, such as BIC, $\pcc$ and $\icc$.
Conversely, transitioning to heavy-tailed distributions results in an increase in the sample noise eigenvalue $\widehat{\lambda}_{r_0+1}$, which impacts estimators that depend on small noise, such as AIC.

\item
In both Sections \ref{subsubsec:study2_high_rank} and \ref{subsubsec:study2_heavy_tailed}, the change in conditions either reduces the sample signal eigenvalue or increases the sample noise eigenvalue.
These changes inevitably affect the methods based on the information criterion, whose selection consistency is determined by the gap conditions.
Despite these challenges, ACT consistently delivers high performance, showcasing its potential superiority under specific conditions.
This not only suggests the potential advantage of ACT but also highlights the limitations of information criterion-based estimators.

 \item
Based on the findings in Sections \ref{subsubsec:study2_sparse_covariance} to \ref{subsubsec:study2_power_spiked}, a noteworthy advantage of information criterion-based estimators lies in their robustness against diverse modifications in the structure of the population covariance matrix $\bm{\Sigma}$.
Whether sparsity is imposed on $\bm{\Sigma}$, an extremely large population eigenvalue $\lambda_1$ is introduced, or a power spiked model is assumed for $\lambda_1, \ldots, \lambda_{r_0}$, the information criterion-based estimators consistently exhibit commendable robustness, unlike several other rank estimation methods within Categories 2 and 3.
\end{itemize}

%
%
%
%
\section{Concluding discussions}\label{sec:conclude}

In this article, we derived the gap conditions for estimators proposed by \citet{bai_ng2002} and further integrated these with the findings of \citet{bai2018} and \citet{hung2022} to develop a unified selection consistency theorem for information criterion-based estimators.
An intriguing observation is that, despite these different methods originating from different working models and distinct formulas, their selection consistency is all governed by the same underlying principle, namely the gap conditions.
This unified selection consistency theorem not only consolidates our understanding of these methods but also offers insightful implications for their application to practitioners and analysts.
We further empirically validated the unified selection consistency theorem through extensive numerical experiments.

Furthermore, we carried out an additional simulation study to explore the characteristics of information criterion-based estimators in comparison with estimators from other categories. 
Specifically, in Sections \ref{subsubsec:study2_high_rank} to \ref{subsubsec:study2_power_spiked}, we observed the impacts on estimation through the five different scenarios. 
As a result, we gained insights into the advantages and potential limitations of information criterion-based estimators. 
A notable advantage of these estimators is their robustness against structural changes in the population covariance, while certain other methods were affected by these changes. 
On the other hand, a potential limitation of information criterion-based estimators is that their selection consistency is determined by gap conditions. 
These conditions outline the requirements for the minimum signal strength and the maximum acceptable noise level. 
Consequently, the information criterion-based estimators can be susceptible to a decrease in signal strength as well as an increase in noise level. 

This study primarily focused on unraveling the properties of information criterion-based estimators; thus, we did not delve deep into methods in other categories. For future research, we propose expanding the scope to include methods from Categories 2 and 3, incorporating both theoretical considerations and extensive numerical experiments. Such endeavors hold the potential to catalyze further advancements in the field of rank estimation methods.

%
%

%
%


%
%

\section*{Declaration of Generative AI and AI-assisted technologies in the writing process}

During the preparation of this work, the authors used ChatGPT for basic grammar and spelling checks, as well as to enhance the readability of the manuscript.
After using this tool, the authors reviewed and edited the content as needed and take full responsibility for the content of the publication.

%
%

\section*{Acknowledgements}
The authors acknowledge financial support from the National Science and Technology Council, Taiwan, under Grant Numbers MOST 110-2118-M-002 -001 -MY3 and MOST 110-2118-M-001 -007 -MY3.

%
%

\bibliographystyle{elsarticle-harv}
\bibliography{reference}

%
%
%
%

\appendix
\renewcommand{\thesection}{\Alph{section}}

%
%

\section{Proof of the unified selection consistency theorem}\label{sec:proof_of_theorem}

\subsection{Lemmas}
%
%
\begin{lem}
\label{lem:lambda_hat_j_vs_psi_lambda_j}
Assume that the general RMT assumptions (C1)-(C3) are satisfied.
Further assume that $\lambda_{r_0}$ goes to infinity as $p \rightarrow \infty$.
Then, we have
\begin{equation}
\frac{\widehat{\lambda}_j}{\psi(\lambda_j)} \convas 1, ~~ \forall j \leq r_0.
\end{equation}
\end{lem}
\begin{proof}
We begin by noting that, for a general distribution $H$, it holds that $\widehat{\lambda}_j/\lambda_j \convas 1$, as demonstrated by \citet[Lemma 2.2]{bai2018}.
Therefore, it suffices to show that $\psi(\lambda)/\lambda \rightarrow 1$ as $\lambda \rightarrow \infty$.
Let us express $\psi(\lambda)/\lambda$ as follows:
\begin{equation}
\frac{\psi(\lambda)}{\lambda} = 1 + c \int \frac{t}{\lambda-t} dH(t).
\end{equation}
By assumption (C3), $H$ is a bounded distribution.
Let $\overline{\lambda}$ be the upper bound of $H$, i.e., $\overline{\lambda} \eqdef {\rm ess~sup}(H)$.
Since we consider the situation, where $\lambda \rightarrow \infty$, we may assume that $\lambda > (k+1) \overline{\lambda}$, where $k \in \mathbb{N}$.
Next, we evaluate the integral term as follows:
\begin{equation}
0 \leq \int \frac{t}{\lambda-t} dH(t) \leq \frac{\overline{\lambda}}{\lambda-\overline{\lambda}} < \frac{1}{k}.
\end{equation}
By taking the limit as $k \rightarrow \infty$, we conclude that $\psi(\lambda)/\lambda \rightarrow 1$ as $\lambda \rightarrow \infty$, which completes the proof.
\end{proof}
%
%
\begin{lem}[Convergence of $\widehat{\sigma}_q^2$]\label{lem:sigma2_hat_convergence}
Assume the general RMT assumptions (C1)-(C3).
For $\widehat{\sigma}_q^2 \eqdef \frac{1}{p-q} \sum_{j=q+1}^p \widehat{\lambda}_j$ with $r_0 \leq q = o(p)$, we have $\widehat{\sigma}_q^2 \convas \mu_H$.
\end{lem}
\begin{proof}
It is sufficient to show the case when $q=r_0$, i.e., $\widehat{\sigma}_q^2 \stackrel{\rm eq}{=} \frac{1}{p-r_0} \sum_{j=r_0+1}^p \widehat{\lambda}_j$.
Define $\widehat{H}(t) \eqdef \frac{1}{p-r_0} \sum_{j=r_0+1}^p \I(\widehat{\lambda}_j \leq t)$.
Since $F^{\bm{S}_n} \convd \F_{c,H}$ with probability 1, we also have $\widehat{H}(t) \convd \F_{c,H}$ with probability 1.
Note that $\widehat{\lambda}_{r_0+1} \convas b$.
Fix $\omega \in \Omega$ s.t $\widehat{\lambda}_{r_0+1} \rightarrow b$.
When $n, p$ are sufficiently large, $\widehat{\lambda}_{r_0+1} \leq b+\epsilon_\omega$ for some $\epsilon_\omega > 0$.
So
\begin{equation}
\widehat{\sigma}_q^2 = \int t d\widehat{H}(t) = \int \min\{t,b+\epsilon_\omega\} d\widehat{H}(t).
\end{equation}
Since $t \mapsto \min\{t, b+\epsilon_\omega\}$ is a bounded continuous function and $\widehat{H}(t) \convd \F_{c,H}$, we have $\widehat{\sigma}_q^2 \rightarrow \int t d\F_{c,H}$.
Furthermore, $\int t d\mathcal{F}_{c,H}(t) = \int t dH(t)$ holds as mentioned in \cite[Lemma 2.16]{yao_zheng_bai_2015}.
Now we have the desired conclusion.
\end{proof}

\subsection{Proof of Theorem~\ref{thm:info_gap_conditions}}
\begin{proof}
The proofs for AIC, BIC and GIC can be found in \citet[Theorems 3-5]{hung2022}.
In the following proof, we discuss the necessary and sufficient conditions (the gap conditions) under which $\widehat{r}_{\rm PC3}$ and $\widehat{r}_{\rm IC3}$ attain strong consistency. 

{\it We first derive the gap conditions for $\widehat{r}_{\rm PC3}$}.
Recall that $\widehat{r}_{\rm PC3}$ is given by $\argmin_{r \leq q} P_3(r)$, where
\begin{equation}
P_3(r) \eqdef \left(\frac{1}{p} \sum_{j=r+1}^p \widehat{\lambda}_j\right) + r \widehat{\sigma}_q^2 \frac{\ln m}{m}, ~~ {\rm and}~ m \eqdef n \wedge p.
\end{equation}
The estimator $\widehat{\sigma}_q^2$ adopted by \citet{bai_ng2002} is given as $\widehat{\sigma}_q^2 \eqdef \frac{1}{p-q} \sum_{j=q+1}^p \widehat{\lambda}_j$, where $q=o(p)$ is the uppper bound of the search range for $r$.
Let $\delta_r \eqdef P_3(r)-P_3(r-1)$, which represents the change in the $\pcc$ criterion as the model rank increases from $r-1$ to $r$.
By rearranging the terms, $\delta_r$ can be expressed as follows:
\begin{eqnarray}
\delta_r &=& \left(\frac{1}{p}\sum_{j=r+1}^p \widehat{\lambda}_j\right) + r \widehat{\sigma}_q^2 \frac{\ln m}{m}
- \left(\frac{1}{p}\sum_{j=r}^p\widehat\lambda_j\right) - (r-1)\widehat{\sigma}_q^2 \frac{\ln m}{m}\\
 &=& -\frac{1}{p} \widehat{\lambda}_r + \widehat{\sigma}_q^2 \frac{\ln m}{m}.
\end{eqnarray}
To ensure that $P_3(r)$ reaches its minimum at $r=r_0$, it is necessary and sufficient to have $\delta_1, \ldots, \delta_{r_0} < 0$ and $\delta_{r_0+1}, \delta_{r_0+2}, \ldots > 0$.
Since $\widehat{\lambda}_r$ is monotonically decreasing, then $\delta_r$ is monotonically increasing, with respect to $r$. Thus, we only need to check that $\delta_{r_0} < 0$ and $\delta_{r_0+1} > 0$.
Hence, the gap conditions for $\pcc$ can be obtained as follows:
\begin{equation}
\begin{cases}
\delta_{r_0} <0 , \\[1ex]
\delta_{r_0 + 1} >0,
\end{cases}
\stackrel{\rm equivalent~to}{\Longleftrightarrow} ~~
\begin{cases}
\frac{p \ln m}{m} <\frac{\widehat{\lambda}_{r_0}}{\widehat{\sigma}_q^2} , \\[1ex]
\frac{p \ln m}{m} > \frac{\widehat{\lambda}_{r_0+1}}{\widehat{\sigma}_q^2}.
\end{cases}
\end{equation}
The first gap condition requires that ${\widehat{\lambda}_{r_0}}/{\widehat{\sigma}_q^2}$ diverges to infinity faster than ${p \ln m}/{m}$.
Using Lemma~\ref{lem:sigma2_hat_convergence}, we have $\widehat{\sigma}_q^2 \convas \mu_H$.
According to Lemma \ref{lem:lambda_hat_j_vs_psi_lambda_j}, we know that $\widehat{\lambda}_{r_0}$ and $\psi_{r_0}$ diverge to infinity at the same order in the scenarios where $\lambda_{r_0} \rightarrow \infty$.
Therefore, we can replace $\widehat{\lambda}_{r_0}$ with $\psi_{r_0}$.
It is important to note that the second gap condition is automatically satisfied.
Since $g$ is a strictly increasing function, the gap conditions for $\pcc$ can be written as follows:
\begin{equation}
\begin{cases}
\frac{p \ln m}{m} <\frac{\psi_{r_0}}{\mu_H} , \\[1ex]
\frac{p \ln m}{m} > \frac{b}{\mu_H},
\end{cases}
\stackrel {\rm equivalent~to}\Longleftrightarrow
~
\begin{cases}
g\left(\frac{\psi_{r_0}}{\mu_H}\right)-g \left(\frac{p \ln m}{m}\right)>0, \\[1ex]
g\left(\frac{b}{\mu_H}\right)-g\left(\frac{p \ln m }{m} \right) < 0.
\end{cases}
\end{equation}
\vspace{13pt}

{\it Next, we derive the gap conditions for $\widehat{r}_{\rm IC3}$}.
Recall that $\widehat{r}_{\rm IC3} \eqdef \argmin_{r \leq q} I_3(r)$, where
\begin{equation}
I_3(r) \eqdef \ln\left(\frac{1}{p} \sum_{j=r+1}^p \widehat{\lambda}_j\right) + r \frac{\ln m}{m}.
\end{equation}
Define $\delta_r \eqdef I_3(r) - I_3(r-1)$.
However, $\delta_r$ is not monotonically increasing with respect to $r$, unlike the $\pcc$ case.
Therefore, we divide the proof into the following two steps:
\begin{description}
  \item[Step 1.]
   We show that the first gap condition is necessary and sufficient for $I_3(r)$ to reach its minimum at $r=r_0$ when $1 \leq r \leq r_0$.
  \item[Step 2.]
   We show that the second gap condition is necessary and sufficient for $I_3(r)$ to reach its minimum at $r=r_0$ for $r_0 \leq r = o(p)$.
\end{description}
%
%
(\textbf{Step~1a}) Here, we derive a necessary condition for the inequality $I_3(r) - I_3(r_0) > 0$ to hold for all $r$ satisfying $1 \leq r \leq r_0-1$.
Furthermore, we will demonstrate that this necessary condition is also sufficient (in \textbf{Step~1b}).
Let us define $\widehat{\mu}_H \eqdef \frac{1}{p} \sum_{j=r_0+1}^p \widehat{\lambda}_j$.
It is worth noting that the expression $I_3(r)-I_3(r_0)$ can be represented as follows:
\begin{equation} \label{eq:I3r_I3r0}
I_3(r) - I_3(r_0) = \ln\left(1 + \frac{\sum_{j=r+1}^{r_0} \widehat{\lambda}_j}{\widehat{\mu}_H p} \right) - (r_0 -r) \frac{\ln m}{m}.
\end{equation}
Again, our objective is to obtain the necessary and sufficient conditions for the validity of $I_3(r) - I_3(r_0) > 0$ for all $r$ satisfying $1 \leq r \leq r_0-1$.
To achieve this, we start by substituting $r = r_0-1$ and derive a necessary condition expressed by the following equation.
(It will be demonstrated later that the condition derived from here is indeed a sufficient condition.)
\begin{equation}\label{eq:I3r0_1_sub_I3r0}
I_3(r_0-1)-I_3(r_0) = \ln\left(1+\frac{\widehat{\lambda}_{r_0}}{\widehat{\mu}_H p}\right)-\frac{\ln m}{m} > 0.
\end{equation}
Equation (\ref{eq:I3r0_1_sub_I3r0}) holds when $\widehat{\lambda}_{r_0} \neq o(p)$ since the term $\frac{\ln m}{m}$ approaches zero,
while $\ln\left(1+\frac{\widehat{\lambda}{r_0}}{\widehat{\mu}_H p}\right) \geq \ln\left(1+\frac{\alpha}{\mu_H}\right)$ eventually holds for some positive number $\alpha$.
Conversely, when $\widehat{\lambda}_{r_0} = o(p)$ (i.e., $\frac{\widehat{\lambda}_{r_0}}{p} \rightarrow 0$), Equation (\ref{eq:I3r0_1_sub_I3r0}) can be approximated as:
\begin{equation}\label{eq:I3r0_1_sub_I3r0_approx}
I_3(r_0-1)-I_3(r_0) \approx \frac{\widehat{\lambda}_{r_0}}{\widehat{\mu}_H p}-\frac{\ln m}{m} > 0,
\end{equation}
because $\ln(1+x) \approx x$ if $|x|$ is sufficiently small.
Hence, we establish that the condition expressed by the inequality below is necessary for $I_3(r)-I_3(r_0) > 0$ to hold for all $r$ satisfying $1 \leq r \leq r_0-1$:
\begin{equation}\label{eq:IC3_first_gap_not_simplified}
\frac{\widehat{\lambda}_{r_0}}{\widehat{\mu}_H} > \frac{p \ln m}{m}.
\end{equation}
(\textbf{Step 1b})
It is further claimed that condition~(\ref{eq:IC3_first_gap_not_simplified}) is sufficient to guarantee $I_3(r)-I_3(r_0) > 0$ for all $r$ with $1 \leq r \leq r_0-1$.
By~(\ref{eq:I3r_I3r0}), we have
\begin{equation}\label{I3r_I3r0_inequality}
I_3(r)-I_3(r_0) \geq \ln\left(1+\frac{(r_0-r)\widehat{\lambda}_{r_0}}{\widehat{\mu}_H p}\right)-(r_0-r)\frac{\ln m}{m}.
\end{equation}
For all $r$ satisfying $1 \leq r \leq r_0-1$, the right-hand side of Equation (\ref{I3r_I3r0_inequality}) is greater than $0$ when $\widehat{\lambda}_{r_0} \neq o(p)$.
This occurs because as ${\ln m}/{m} \rightarrow 0$, $\ln\left(1+\frac{(r_0-r)\widehat{\lambda}_{r_0}}{\widehat{\mu}_H p}\right)$ eventually becomes positive.
On the other hand, when $\widehat{\lambda}_{r_0} \neq o(p)$, the right-hand side of the inequality can be approximated as follows:
\begin{equation}
\frac{(r_0-r)\widehat{\lambda}_{r_0}}{\widehat{\mu}_H p} - (r_0-r)\frac{\ln m}{m} = \frac{(r_0-r)}{p} \left(\frac{\widehat{\lambda}_{r_0}}{\widehat{\mu}_H} - \frac{p \ln m}{m}\right).
\end{equation}
This expression is positive when Equation (\ref{eq:IC3_first_gap_not_simplified}) holds.
Therefore, we can conclude that Equation (\ref{eq:IC3_first_gap_not_simplified}) serves as both a necessary and sufficient condition for $I_3(r)-I_3(r_0)$ to be positive for all $r$ with $1 \leq r \leq r_0-1$.
Similar to the proof of $\pcc$, this condition can be equivalently represented by the following equation:
\[
g\left(\frac{\psi_{r_0}}{\mu_H}\right) - g\left(\frac{p \ln m}{m}\right) > 0.
\]
%
%
(\textbf{Step2})  We derive a necessary and sufficient condition $I_3(r)-I_3(r_0) > 0$ for all $r > r_0$ with $r=o(p)$.
When $r_0 < r = o(p)$, we have
\begin{eqnarray*}
I_3(r) - I_3(r_0) &=& \ln\left(\frac{1}{p} \sum_{j=r+1}^p \widehat{\lambda}_j\right) - \ln\left(\frac{1}{p} \sum_{j=r_0+1}^p \widehat{\lambda}_j\right) + (r-r_0) \frac{\ln m}{m} \\
&=& \ln\left(1 - \frac{\sum_{j=r_0+1}^r \widehat{\lambda}_j}{p \widehat{\mu}_H}\right) + (r-r_0) \frac{\ln m}{m} \\
&\approx& \ln\left(1 - \frac{(r-r_0)b}{\mu_H p}\right) + (r-r_0) \frac{\ln m}{m} \\
&\approx&  - \frac{(r-r_0)b}{p \mu_H} + (r-r_0) \frac{\ln m}{m},
\end{eqnarray*}
where the first approximation is based on the fact that $\widehat{\lambda}_r \rightarrow b$ a.s., when $r_0 < r=o(p)$ and
the second approximation is based on $\ln(1-x) \approx -x$ for $|x| =o(1)$.
We obtain $I_3(r)-I_3(r_0) > 0$ if and only if
\[
\frac{p \ln m}{m}> \frac{b}{\mu_H} ~ \stackrel{\rm equivalent~to}\Leftrightarrow ~ g\left(\frac{b}{\mu_H}\right)-g\left(\frac{p \ln m}{m}\right)<0.
\]
Certainly, it is true that ${p \ln m}/{m}$ diverges to infinity, so the inequality mentioned above always holds.
By performing the transformations as described above, we are able to express the gap conditions for all information-criterion estimators, including AIC, BIC, and GIC, in a unified formulation as in Equation (\ref{eq:unified_gap_conditions}).
(Note that as $n, p \rightarrow \infty$, we have ${p}/{m} \rightarrow (1 \vee c)$.)

For $\pca$/$\ica$ and $\pcb$/$\icb$, we can derive the following gap conditions by following a similar line of reasoning as discussed above.
Moreover, these can also be represented in the form of Equation (\ref{eq:unified_gap_conditions}) through simple transformations.
\[
(\pca/\ica)~
\begin{cases}
\frac{\psi_{r_0}}{\mu_H} > \left(1+\frac{p}{n}\right) \ln\left(\frac{np}{n+p}\right), \\[1ex]
\frac{b}{\mu_H} <  \left(1+\frac{p}{n}\right) \ln\left(\frac{np}{n+p}\right) ,
\end{cases}
~~~ (\pcb/\icb)~
\begin{cases}
\frac{\psi_{r_0}}{\mu_H} > \left(1+\frac{p}{n}\right) \ln(m), \\[1ex]
\frac{b}{\mu_H} <  \left(1+\frac{p}{n}\right) \ln(m).
\end{cases}
\]
Their proofs are similar to those for $\pcc$ and $\icc$ and thus omitted.
\end{proof}

%
%

\section{Study III:  Comparison of BIC and PCs/ICs in the regime of $p\gg n$}\label{sec:study3}

In Section \ref{subsec:study1}, we conducted simulations for $c=0.25, 0.5, 1.0, 1.5$.
We believe that Section \ref{subsec:study1} offers a sufficiently extensive set of settings and their corresponding results to numerically validate the unified selection consistency theorem.
However, we also provide simulation results for larger values of $c$ in Appendix \ref{sec:study3}.
Specifically, in this section, we will focus on the performance of BIC, PCs, and ICs.

\subsection{Experimental settings}

The simulation procedures used in this section closely align with those implemented in Section \ref{subsubsec:study1_settings}.
Consistent with Section \ref{subsubsec:study1_settings}, we keep the value of $r_0$ steady at 5.
However, we examine different combinations of $(n,p)$, specifically $(300, 900)$, $(250, 1000)$, and $(240, 1200)$.
The limiting spectral distribution (LSD) of the population covariance matrix $\bm{\Sigma}$ is represented by $H_4(t) \eqdef \mathcal{I}(1 \leq t)$.
Furthermore, we set $\lambda_{r_0}$ to various values: 5, 10, 15, 20, 25, 30, and 35.

\subsection{Gap conditions}

\begin{figure}[h!]
\includegraphics[width=12cm]{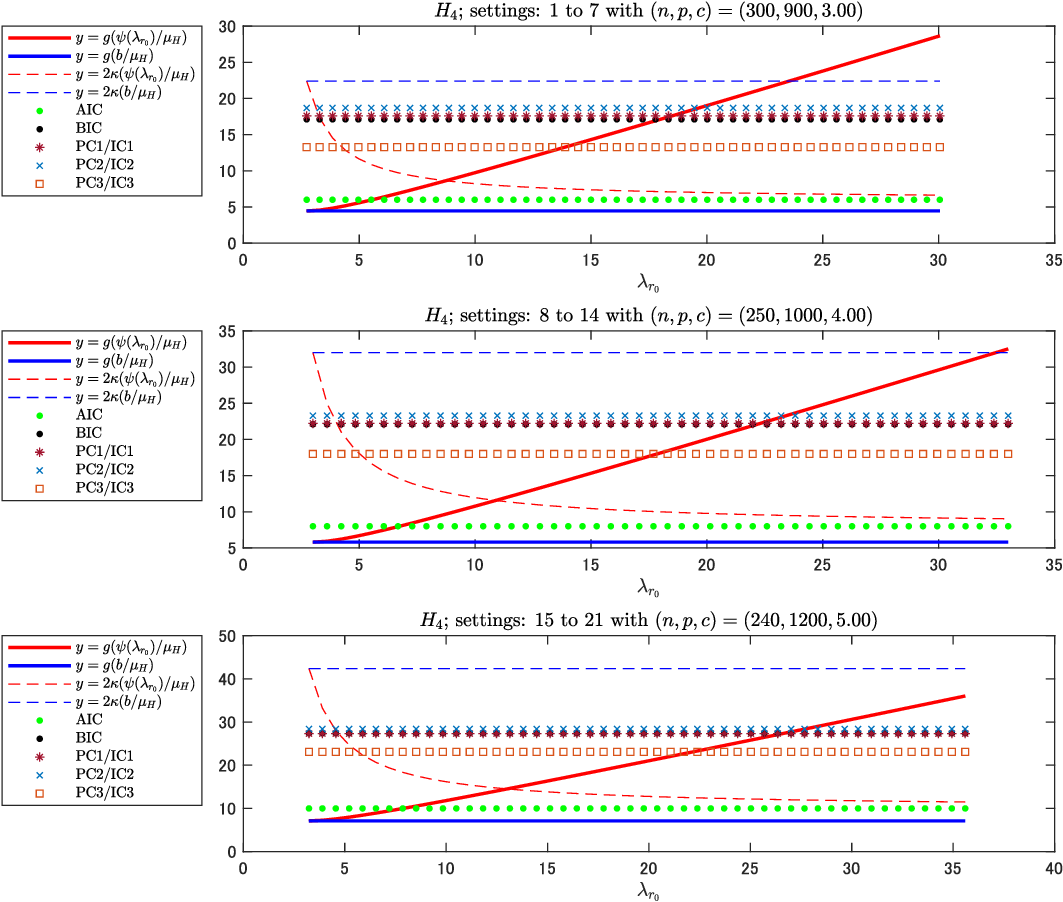}
\centering
\caption{Gap condition graphs for $p\gg n$.
The three graphs correspond to $(n,p)$ values of  (300,900), (250,1000), and (240,1200), respectively.
Additionally, they correspond to settings 1-7, 8-14, and 15-21 in Table \ref{table:large_p_case_simulation}, respectively.
The gap condition value for $\pcb/\icb$ is the largest, followed by BIC and $\pca/\ica$, whose gap condition values are quite close.
On the other hand, the gap condition value for $\pcc/\icc$ is relatively small compared to those of BIC, PCs/ICs.
}
\label{fig:graph_of_gap_conditions_large_p_case}
\end{figure}

Figure \ref{fig:graph_of_gap_conditions_large_p_case} displays the gap condition graphs for the aformentioned settings.
The individual graphs correspond to the settings where the pairs $(n,p)$ are set as $(300, 900)$, $(250, 1000)$, and $(240, 1200)$, respectively.
As elaborated in Theorem \ref{thm:info_gap_conditions}, we examined the interrelationships among the $\beta_{\rm m}$ functions.
For $c > 1$, the inequality $\beta_3 < \beta_{\rm BIC} < \beta_1 < \beta_2$ is satisfied.
Additionally, as $c$ rises, $\beta_1$ narrows its gap with $\beta_2$.
On the contrary, in our chosen settings, specifically where $n=250$ or $n=240$, $\beta_{\rm BIC}$ is close to $\beta_1$, and even exceeds it, due to the insufficiently large value of $n$.
Yet, for sufficiently large $n$, the relationship $\beta_{\rm BIC} < \beta_1$ is expected to hold.
Drawing upon these observations, we expect that $\pcc/\icc$ can achieve high accuracy with relatively small signal eigenvalues, while BIC, $\pcb/\icb$, and $\pcc/\icc$ require stronger signal eigenvalues than $\pcc/\icc$.
Additionally, the performance of BIC, $\pcb/\icb$, and $\pcc/\icc$ is expected to be quite similar

\subsection{Results and comments}\label{subsec:large_p_results}

The simulation results are depicted in Table \ref{table:large_p_case_simulation}.
As in previous simulations, we do not estimate $\sigma^2$ for $\pca, \pcb$, and $\pcc$.
Instead, we directly substitute $\widehat{\sigma}_q^2 \leftarrow \mu_H \stackrel{\rm eq}{=} 1$.
As evidenced by the results in Table \ref{table:large_p_case_simulation}, $\pcc/\icc$ require the smallest values of $\lambda_{r_0}$ to achieve high estimation accuracy, due to having the lowest $\beta_{\rm m}$ function values (excluding AIC and GIC).
Conversely, $\pcb/\icb$ require larger values of $\lambda_{r_0}$ to achieve similar accuracy, as they have the highest $\beta_{\rm m}$ function values.
While $\pca$ and $\ica$, as well as $\pcb$ and $\icb$ and $\pcc$ and $\icc$, share the same gap conditions, slight performance differences exist between them.
(As discussed in Section \ref{subsubsec:study2_high_rank}, $\pca$ imposes a larger penalty term than $\ica$, resulting in this difference.)
In our current settings, the values of the $\beta_{\rm m}$ functions for BIC and $\pca/\ica$ are quite close, showing similar performance overall.

\begin{table}[h!]
\scalebox{0.8}{
\begin{tabular}{|c|c|c|c|c|c|c|c|c|c|c|c|c|c|}
\hline
   & $n$ & $p$  & $c$ & $\lambda_{r_0}$ & AIC   & BIC   & GIC   & $\pca$ & $\ica$ & $\pcb$ & $\icb$ & $\pcc$ & $\icc$ \\ \hline
1  & 300 & 900  & 3   & 5               & 0.000 & 0.000 & 0.000 & 0.000  & 0.000  & 0.000  & 0.000  & 0.000  & 0.000  \\ \hline
2  & 300 & 900  & 3   & 10              & 1.000 & 0.000 & 0.892 & 0.000  & 0.000  & 0.000  & 0.000  & 0.000  & 0.000  \\ \hline
3  & 300 & 900  & 3   & 15              & 1.000 & 0.000 & 1.000 & 0.000  & 0.000  & 0.000  & 0.000  & 0.038  & 0.084  \\ \hline
4  & 300 & 900  & 3   & 20              & 1.000 & 0.126 & 1.000 & 0.034  & 0.054  & 0.000  & 0.002  & 0.998  & 1.000  \\ \hline
5  & 300 & 900  & 3   & 25              & 1.000 & 0.988 & 1.000 & 0.956  & 0.974  & 0.764  & 0.854  & 1.000  & 1.000  \\ \hline
6  & 300 & 900  & 3   & 30              & 1.000 & 1.000 & 1.000 & 1.000  & 1.000  & 1.000  & 1.000  & 1.000  & 1.000  \\ \hline
7  & 300 & 900  & 3   & 35              & 1.000 & 1.000 & 1.000 & 1.000  & 1.000  & 1.000  & 1.000  & 1.000  & 1.000  \\ \hline
8  & 250 & 1000 & 4   & 5               & 0.000 & 0.000 & 0.000 & 0.000  & 0.000  & 0.000  & 0.000  & 0.000  & 0.000  \\ \hline
9  & 250 & 1000 & 4   & 10              & 0.998 & 0.000 & 0.032 & 0.000  & 0.000  & 0.000  & 0.000  & 0.000  & 0.000  \\ \hline
10 & 250 & 1000 & 4   & 15              & 1.000 & 0.000 & 0.994 & 0.000  & 0.000  & 0.000  & 0.000  & 0.000  & 0.000  \\ \hline
11 & 250 & 1000 & 4   & 20              & 1.000 & 0.000 & 1.000 & 0.000  & 0.000  & 0.000  & 0.000  & 0.064  & 0.146  \\ \hline
12 & 250 & 1000 & 4   & 25              & 1.000 & 0.152 & 1.000 & 0.072  & 0.140  & 0.012  & 0.028  & 0.956  & 0.978  \\ \hline
13 & 250 & 1000 & 4   & 30              & 1.000 & 0.950 & 1.000 & 0.904  & 0.934  & 0.682  & 0.790  & 1.000  & 1.000  \\ \hline
14 & 250 & 1000 & 4   & 35              & 1.000 & 0.998 & 1.000 & 0.998  & 0.998  & 0.988  & 0.990  & 1.000  & 1.000  \\ \hline
15 & 240 & 1200 & 5   & 5               & 0.000 & 0.000 & 0.000 & 0.000  & 0.000  & 0.000  & 0.000  & 0.000  & 0.000  \\ \hline
16 & 240 & 1200 & 5   & 10              & 0.844 & 0.000 & 0.000 & 0.000  & 0.000  & 0.000  & 0.000  & 0.000  & 0.000  \\ \hline
17 & 240 & 1200 & 5   & 15              & 1.000 & 0.000 & 0.858 & 0.000  & 0.000  & 0.000  & 0.000  & 0.000  & 0.000  \\ \hline
18 & 240 & 1200 & 5   & 20              & 1.000 & 0.000 & 1.000 & 0.000  & 0.000  & 0.000  & 0.000  & 0.000  & 0.000  \\ \hline
19 & 240 & 1200 & 5   & 25              & 1.000 & 0.000 & 1.000 & 0.000  & 0.000  & 0.000  & 0.000  & 0.060  & 0.130  \\ \hline
20 & 240 & 1200 & 5   & 30              & 1.000 & 0.092 & 1.000 & 0.062  & 0.108  & 0.022  & 0.034  & 0.858  & 0.912  \\ \hline
21 & 240 & 1200 & 5   & 35              & 1.000 & 0.856 & 1.000 & 0.800  & 0.872  & 0.562  & 0.692  & 1.000  & 1.000  \\ \hline
\end{tabular}
}
\centering
\caption{Simulation results in the regime of $p \gg n$.
As depicted in Figure \ref{fig:graph_of_gap_conditions_large_p_case}, the value of the $\beta_{\rm m}(\cdot)$ function for $\pcb/\icb$ is the highest.
It is followed closely by BIC and $\pca/\ica$, which exhibit nearly identical values, while $\pcc/\icc$ shows relatively lower values.
The simulation results corroborate these observations. Specifically, $\pcc/\icc$ attains selection consistency even with relatively weak signals.
Meanwhile, BIC and $\pca/\ica$ demonstrate similar performance throughout the simulation.
Furthermore, it is evident that $\pcb/\icb$ necessitates the strongest signals for comparable performance.
}
\label{table:large_p_case_simulation}
\end{table}


%
%

\section{Some selected rank estimation methods from Categories 2 and 3}\label{sec:intro_category2_and_3}


\subsection{Category 2: ACT, DPA and BEMA}

The methods within this category share a common characteristic of estimating the essential supremum of the LSD of the sample covariance matrix $\bm{S}_n$ and utilizing it for rank estimation.
Alongside ACT \citep{act} and DPA \citep{dpa}, another method known as BEMA \citep{bema} is also included.

Despite ACT employing the sample correlation matrix $\widehat{\bm{R}}$ instead of the sample covariance matrix $\bm{S}_n$, we include it within this category.
Let $\{\widehat{\lambda}_j\}_{j=1}^p$ be the eigenvalues of the sample correlation matrix $\widehat{\bm{R}}$.
The estimator by ACT is expressed by Equation (\ref{eq:act_estimator}).
\begin{equation}\label{eq:act_estimator}
\widehat r_{\rm ACT} \eqdef \max\{j \mid \widetilde\lambda_j>1+\sqrt{\widehat{c}}\},
\end{equation}
where $\tilde{\lambda}_j$ is given by the following equation:
\begin{eqnarray}
\tilde{\lambda}_j &\eqdef& -\frac{1}{\underline{m}_{n,j}(\widehat{\lambda}_j)}, \\
\underline{m}_{n,j}(z) &\eqdef& -\left(1-\frac{p-j}{n-1}\right) z^{-1} + \left(\frac{p-j}{n-1}\right) m_{n,j}(z), \\
m_{n,j}(z) &\eqdef& \frac{1}{p-j}\left\{
\sum_{\ell=j+1}^p (\widehat{\lambda}_\ell-z)^{-1} + \left(\frac{3 \widehat{\lambda}_j}{4}+\frac{\widehat{\lambda}_{j+1}}{4} - z \right)^{-1} \right\}.
\end{eqnarray}

DPA \citep{dpa} is known as the deterministic version of Parallel Analysis \citep{horn1965}.
It estimates the essential supremum of the LSD of $\bm{S}_n$ using Equation (\ref{eq:dpa_b_hat}):
\begin{equation}\label{eq:dpa_b_hat}
\widehat{b}_{\rm DPA} \eqdef {\rm ess~sup}(\F_{p/n,\widehat{H}}),
\end{equation}
where $\widehat{H}$ is defined by the empirical spectral distribution (ESD) of ${\rm diag}\left(\bm{S}_n\right)$. The rank is then estimated using Equation (\ref{eq:dpa_estimator}):
\begin{equation}\label{eq:dpa_estimator}
\widehat{r}_{\rm DPA} \eqdef \sum_{j=1}^{p}\mathcal{I}\left(\widehat\lambda_j>\widehat{b}_{\rm DPA}\right),
\end{equation}
where $\{\widehat{\lambda}_j\}_{j=1}^p$ represents the eigenvalues of $\bm{S}_n$.

BEMA \citep{bema} is based on a general spiked model represented by Equation (\ref{eq:bema_model}).
It assumes a Gamma distribution for the variances of the noise, which allows for non-constant variances.
\begin{equation}\label{eq:bema_model}
\begin{cases}
\bm{\Sigma} \eqdef \sum_{j=1}^r \mu_j \bm{\gamma}_j\bm{\gamma}_j^\top + {\rm diag}(\sigma_1^2,\ldots,\sigma_p^2), \\
{\rm where} \quad \sigma_1^2,\ldots,\sigma_p^2 \iid {\rm Gamma}(\theta,\theta/\sigma^2).
\end{cases}
\end{equation}
In this model, the LSD of $\bm{\Sigma}$ is also $H \eqdef {\rm Gamma}(\theta,\theta/\sigma^2)$ as it eventually can ignore the influence of the signal eigenvalues.
Let $\F_{p/n,H}$ represent the generalized Marchenko-Pastur distribution, which corresponds to the LSD of $\bm{S}_n$.
The rank estimation method of BEMA involves two steps.
\begin{itemize}\itemsep=0pt
  \item First, the bulk sample eigenvalues $\{\widehat{\lambda}_j\}_{j \in [\alpha (n \wedge p), (1-\alpha)n \wedge p]}$ are fitted to $\F_{p/n, H}$ to estimate the unknown parameters $(\sigma^2, \theta)$.
  \item Second, the essential supremum of the distribution $\F_{p/n, H}$ is estimated using Monte Carlo estimation.
\end{itemize}
The following procedure is repeated $M$ times, with $m$ denoting the current iteration ($1 \leq m \leq M$):
\begin{enumerate}
\item Generate a noise variance matrix $\bm{D}_m \eqdef {\rm diag}(d^{(m)}_1, \ldots, d^{(m)}_p)$ based on the estimated Gamma distribution,
i.e., $d^{(m)}_1,\ldots,d^{(m)}_p \iid {\rm Gamma}(\widehat{\theta}, \widehat{\theta}/\widehat{\sigma}^2)$.
\item Generate a sample of size $n$ from a multivariate normal distribution with the generated noise variance $\bm{D}_m$, i.e., $\bm{x}^{(m)}_1, \ldots, \bm{x}^{(m)}_n \iid N(\bm{0}, \bm{D}_m)$.
\item Find the largest eigenvalue $\widehat{\lambda}_1(m)$ of the sample covariance matrix $\frac{1}{n} \bm{X}_m^\top \bm{X}_m$, where $\bm{X}_m\eqdef [\bm{x}^{(m)}_1,\ldots,\bm{x}^{(m)}_n]^\top$.
\item Update $m \leftarrow m+1$ and repeat steps 1-3 unless $m = M$.
\end{enumerate}
Finally, let $\widehat{b}$ be the $(1-\beta)$th quantile of $\{\widehat{\lambda}_1(m)\}_{m=1}^M$ and estimate $r_0$ by
\begin{equation}
\widehat r_{\rm BEMA} = \sum_{j=1}^{p}\mathcal{I}\{\widehat\lambda_j> \widehat{b}\}.
\end{equation}

\subsection{Category 3: ED, ON and ER/GR}

If $1 \leq j \leq r_0$, there exists a noticeable gap between $\widehat{\lambda}_j$ and $\widehat{\lambda}_{j+1}$.
On the other hand, if $r_0 < j \leq q$, both $\widehat{\lambda}_j$ and $\widehat{\lambda}_{j+1}$ converge to $b$.
This observation underscores the effectiveness of examining the difference between adjacent sample eigenvalues in estimating the rank.

ED \citep{ed} focuses on the difference between neighboring sample eigenvalues to estimate $r_0$.
The estimation is given by Equation (\ref{eq:ed_estimator}):
\begin{eqnarray}\label{eq:ed_estimator}
\widehat r_{\rm ED} \eqdef \max\{j\le q: \widehat\lambda_j-\widehat\lambda_{j+1}\ge \delta\},
\end{eqnarray}
where $\delta$ is a threshold calibrated through the following procedure.
The procedure consists of the following steps, which are repeated until convergence:
(Initially, $j$ is set to $q+1$, where $q$ is the upper bound of the search range for $r$)
\begin{enumerate}
 \item Perform a simple linear regression from ${(j-1)^{2/3},\ldots,(j+3)^{2/3}}$ to ${\lambda_j, \ldots, \lambda_{j+4}}$.
 Obtain the slope coefficient $\widehat{\beta}$ using the least squares method and define $\delta \eqdef 2|\widehat{\beta}|$.
 \item Compute $\widehat{r} \eqdef \max\{j \leq q \mid \widehat{\lambda}_j - \widehat{\lambda}_{j+1} \geq \delta\}$ where $\max\emptyset$ is defined to be 0.
 \item Update $j \leftarrow \widehat{r}+1$.
\end{enumerate}

The estimators of ON \citep{on}, ER, and GR \citep{er_gr} are defined as follows in Equations (\ref{eq:on_estimator}) through (\ref{eq:gr_estimator}).
\begin{eqnarray}
\widehat r_{\rm ON}&\eqdef&\argmax_{j\le q}\frac{\widehat\lambda_j-\widehat\lambda_{j+1}}{\widehat\lambda_{j+1}-\widehat\lambda_{j+2}}, \label{eq:on_estimator} \\
\widehat r_{\rm ER}&\eqdef&\argmax_{j\le q}\frac{\widehat\lambda_j}{\widehat\lambda_{j+1}}, \label{eq:er_estimator} \\
\widehat r_{\rm GR}&\eqdef&\argmax_{j\le q}\frac{\ln\left(1+\frac{\widehat\lambda_j}{\sum_{\ell>j}\widehat\lambda_\ell}\right)}
{\ln\left(1+\frac{\widehat\lambda_{j+1}}{\sum_{\ell>j+1}\widehat\lambda_\ell}\right)} \label{eq:gr_estimator}.
\end{eqnarray}

\end{document}